\documentclass[12pt]{amsart}

\def\@typesizes{
       \or{5}{6.5}\or{6}{7.5}\or{7}{8.5}\or{8}{11}\or{9}{12}
       \or{10}{13}
       \or{\@xipt}{14}\or{\@xiipt}{15}\or{\@xivpt}{18}
       \or{\@xviipt}{20}\or{\@xxpt}{24}}

\usepackage[utf8]{inputenc}
\usepackage[T1]{fontenc}
\usepackage[english]{babel}
\usepackage{hyperref}
\usepackage{graphicx}
\usepackage{amsmath,amsfonts,amssymb,amsthm,cite}

\numberwithin{equation}{section}
\newtheorem{theorem}{Theorem}[section]
\newtheorem{corollary}[theorem]{Corollary}

\newtheorem{lemma}[theorem]{Lemma}
\newtheorem{remark}[theorem]{Remark}

\newtheorem{example}[theorem]{Example}
\newtheorem{definition}[theorem]{Definition}
\newcommand{\nonprint}[1]{}

\begin{document}

\title[Parameter-dependent boundary-value problems]{Parameter-dependent inhomogeneous boundary-value problems in spaces of continuously differentiable functions}

\author{O. M. Atlasiuk}
\address{
University of Helsinki, Department of Mathematics and Statistics, P.O. Box 68, Pietari Kalmin katu 5, 00014 Helsinki, Finland and \\
Institute of Mathematics of the National Academy of Sciences of Ukraine, st. Tereschenkivska 3, 01024 Kyiv, Ukraine }
\email{olena.atlasiuk@helsinki.fi}

\author{V. A. Mikhailets}
\address{King's College London, Strand, WC2R 2LS London, UK and \\
Institute of Mathematics of the National Academy of Sciences of Ukraine, st. Tereschenkivska 3, 01024 Kyiv, Ukraine }
\email{vladimir.mikhailets@gmail.com}

\begin{abstract}
We investigate a broad class of linear inhomogeneous boundary-value problems of order $r$ for systems of ordinary differential equations depending on a parameter $\mu$ belonging to an abstract metric space $\mathcal M$. The solutions belong to the spaces of continuously differentiable vector-valued functions $(C^{(n+r)})^m$, $n\in\mathbb{N}\cup\{0\}$, $m, r \in \mathbb{N}$. The boundary conditions have the most general form $By=c$, where $B$ is an arbitrary continuous linear operator from $(C^{(n+r)})^m$ to $\mathbb{C}^{rm}$. Thus, they may contain derivatives of the unknown vector-valued function of integer and/or fractional orders up to $n+r$ inclusive. We establish necessary and sufficient conditions for the continuous dependence of solutions on the parameter $\mu$. Furthermore, we prove that the solutions of the original problems can be approximated in the space $(C^{(n+r)})^m$ by solutions of systems of ordinary differential equations with polynomial coefficients, polynomial right-hand sides, and multipoint boundary conditions, where the approximating problems are constructed independently of the right-hand sides of the original problem.
\end{abstract}

\maketitle



\section{Introduction}\label{Sec.1}

Linear boundary-value problems for systems of ordinary differential equations depending on parameters occupy a fundamental place in modern mathematical analysis, control theory, and mathematical physics. The study of the continuous dependence of solutions on parameters is of critical importance for the fidelity of mathematical modeling of numerous real-world processes. Physical parameters of the medium, frequency characteristics, and material properties are invariably determined with a certain error or are described by asymptotic or regularized models. In this connection, a fundamental problem arises: \emph{under what analytical conditions does a perturbation of the parameters of the problem guarantee a small deviation of its solution in the corresponding functional norm?}

Historically, the systematic study of the dependence of solutions of differential equations on a discrete or continuous numerical parameter began for the initial-value (Cauchy) problem in the classical works of I.~I.~Gikhman~\cite{Gikhman}, M.~A.~Krasnoselskii and S.~G.~Krein~\cite{KrasnKrein}, and J.~Kurzweil and Z.~Vorel~\cite{KurzweilVorel}. For linear systems, these results were strengthened and generalized by A.~Yu.~Levin~\cite{Levin}, Z.~Opial~\cite{Opial}, W.~T.~Reid~\cite{Reid}, and T.~K.~Nguyen~\cite{Nguen}.

For boundary-value problems, the investigation of continuity with respect to a parameter has turned out to be fundamentally more intricate due to the non-local nature of boundary conditions and the possible loss of unique solvability. Questions related to limit transitions in general inhomogeneous linear boundary-value problems were thoroughly investigated in the papers of I.~T.~Kiguradze~\cite{Puza, Kigyradze1987, Kigyradze2003, Kigyradze1975} and further developed by his followers~\cite{Ashordia97, Ashordia98, Ashordia, Ashordia05, Ashordia20, Boyko, KodliukMikhailets2013, Pelekhata18}.

For a long time, limit theorems were established exclusively for solutions of boundary-value problems of the form
\begin{equation*}
    Ly = f, \quad By = c,
\end{equation*}
where $L$ is a system of $m$ linear ordinary differential equations of order $r \in \mathbb{N}$, and $B$ is a continuous linear operator on the Banach space of $(r-1)$-times continuously differentiable functions.

In this setting, general boundary conditions can contain only derivatives of the unknown function of orders not exceeding $r-1$~\cite{Krall, Iha, Chua, Paukstaite, Goodrich, Whyburn, Smogorzhevsky, Krein82, Skubachevskii, Ma, Adomian, Krall68}. Sufficient conditions for the continuity of solutions with respect to the parameter were improved and generalized to systems of differential equations of order $r \ge 2$ in the works~\cite{Pelekhata19, KodliukMikhailetsReva2013}. For such problems, the boundary operator $B$ is compact relative to the differential operator $L$, and the problem admits a reduction to classical Fredholm integral equations of the second kind, while the convergence of solutions can be investigated by standard perturbation methods (see, e.g., \cite{Kato} and references therein).

At the same time, in many modern applied problems, boundary conditions naturally arise in which the orders of derivatives \emph{reach or exceed} the order of the differential equation. These include problems in the mechanics of elastic systems and vibration theory, automatic control theory, optimization, and inverse problems.

In such cases, the classical theory becomes fundamentally inapplicable: the operator $B$ is no longer relatively compact, the concept of a Green's function in the classical sense loses its meaning, and the formally adjoint problem as well as the classical Lagrange identity are not defined. This necessitated the development of a new functional-analytic approach, in which the boundary operator is treated as an abstract bounded linear operator
\begin{equation*}
    B \colon (C^{(n+r)})^m \to \mathbb{C}^{rm}, \quad n \in \mathbb{N} \cup \{0\},
\end{equation*}
defined on the solution space of the differential system. Such boundary conditions have been termed generic boundary conditions. This new direction of research was developed in a series of works by one of the authors and his collaborators~\cite{HRYPMikhailetsMurach2017, HRYPKodliukMikhailets2015, HRYP, AtlasiukMikhailets20192, AtlasiukMikhailetsSkor}. The obtained applications included multipoint boundary-value problems~\cite{Ashordia98, Ashordia05, Atlasiuk2021}, as well as problems in the spectral theory of differential operators with highly singular coefficients~\cite{Horyunov}. Subsequently, the class of problems under investigation was significantly broadened by incorporating various function spaces.

However, for systems of \textit{arbitrary order} $r \in \mathbb{N}$ in the spaces of continuously differentiable vector-valued functions $(C^{(n+r)})^m$, depending on a parameter $\mu \in \mathcal{M}$ from an abstract metric space, the following key problems remained open:
\begin{enumerate}
    \item[\emph{(I)}] \emph{Criterion of continuous dependence on an abstract parameter:} To find sharp necessary and sufficient conditions on the coefficients of the system and the boundary operator under which the solution $y(\cdot, \mu)$ depends continuously on the parameter $\mu$ in the space $(C^{(n+r)})^m$.

    \item[\emph{(II)}] \emph{Analytical structure of the abstract boundary operator:} To obtain an explicit canonical Riemann--Stieltjes integral representation for an arbitrary linear operator $B \colon (C^{(n+r)})^m \to \mathbb{C}^{rm}$ and to establish constructive criteria for its strong and uniform operator convergence.

    \item[\emph{(III)}]\emph{ Ulam--Hyers stability:} To establish sharp two-sided estimates for the deviation between the perturbed and unperturbed solutions in terms of the residual of the problem in the context of the modern theory of stability of operator equations (see~\cite{BrzdekPopaRus2018, PopaRasa2011, Jung2011, Rus}).

    \item[\emph{(IV)}] \emph{Constructive approximation by elementary problems:} To clarify the possibility of approximating an arbitrary non-classical boundary-value problem by solutions of systems with polynomial coefficients and multipoint boundary conditions, and to determine the fundamental boundary between strong and uniform operator approximation.
\end{enumerate}

The aim of the present paper is to provide an exhaustive resolution of these formulated problems for linear differential systems of order $r \in \mathbb{N}$ in the spaces $(C^{(n+r)})^m$.

The article is organized as follows. Section~\ref{Sec.2} contains the statements and proofs of general results on the continuous parameter dependence of solutions. In Section~\ref{Sec.3}, we prove limit theorems for the operators appearing in the boundary-value problems. In addition, Theorem~\ref{opBinSobolev} establishes a canonical representation for the abstract boundary operator $B$. This leads to a more detailed and constructive formulation of the results of Section~\ref{Sec.2}. Section~\ref{Sec.4} is devoted to applications of the main results, specifically, theorems on the approximation of solutions of an arbitrary boundary-value problem by solutions of elementary multipoint boundary-value problems.

\section{General results on the continuous parameter dependence}\label{Sec.2}

Let $(a,b)\subset\mathbb{R}$ be a finite open interval, and let the parameters
$$
n\in\mathbb{N}\cup\{0\},\quad\{m, r\}\subset\mathbb{N}
$$
be arbitrary. We denote by
$$C^{(n+r)}:=C^{(n+r)}\bigl([a,b];\mathbb{C}\bigr)$$
the space of $(n+r)$-times continuously differentiable complex-valued functions $y \colon [a,b] \rightarrow \mathbb{C}$. This space is a Banach space equipped with the norm
$$
\big\|y\big\|_{(n+r)}:=\sum_{s=0}^{n+r}\max \big \{\big|y^{(s)}(t) \big| \colon t \in [a,b] \big\}.
$$
We also introduce the spaces
$$
\big(C^{(n+r)}\big)^m:=C^{(n+r)}\bigl([a,b];\mathbb{C}^{m}\bigr)
\;\;\mbox{and}\;\;
\big(C^{(n+r)}\big)^{m\times m}:=C^{(n+r)}\bigl([a,b];\mathbb{C}^{m\times m}\bigr),
$$
which consist, respectively, of vector-valued functions of dimension $m$ and matrix-valued functions of dimension $m\times m$ whose elements belong to the space $C^{(n+r)}$.
The norm in $\big(C^{(n+r)}\big)^m$ is defined as the sum of the norms of its components in $C^{(n+r)}$, and the norm in $\big(C^{(n+r)}\big)^{m\times m}$ is defined as the sum of the norms of its columns in $\big(C^{(n+r)}\big)^m$. We denote these norms by the common symbol $\|\cdot\|_{(n+r)}$, as their precise meaning will always be clear from the context.
The same convention will be adopted for the other function spaces under consideration.

Let $\mathcal{M}$ be an arbitrary metric space. In the sequel, $\mu \in \mathcal{M}$ denotes a variable parameter, whereas $\mu_0 \in \mathcal{M}$ denotes an arbitrarily fixed value of this parameter.
We consider the following linear boundary-value problem for the unknown vector-valued function $y(\cdot,\mu)\in \big(C^{(n+r)}\big)^m$:
\begin{equation}\label{bound_z1}
\left(L(\mu)y(\mu)\right)(t):=y^{(r)}(t,\mu) + \sum\limits_{\ell=1}^rA_{r-\ell}(t,\mu)y^{(r-\ell)}(t,\mu)=f(t,\mu), \quad t\in(a,b),
\end{equation}
\begin{equation} \label{bound_z2}
B(\mu)y(\cdot,\mu)=c(\mu),
\end{equation}
where the matrix-valued functions $A_{r-\ell}(\cdot,\mu) \in \big(C^{(n)}\big)^{m\times m}$, the vector-valued function $f(\cdot,\mu) \in \big(C^{(n)}\big)^m$, the vector $c(\mu) \in \mathbb{C}^{rm}$, and the continuous linear operator
\begin{equation}\label{oper_B(e)v}
B(\mu)\colon (C^{(n+r)})^m\rightarrow\mathbb{C}^{rm}
\end{equation}
are given and arbitrary.
The boundary condition \eqref{bound_z2} consists of $rm$ scalar conditions for a system of $m$ differential equations of $r$-th order. In what follows, we use standard matrix algebra rules; in particular, vectors and vector-valued functions are regarded as column vectors.

A solution of problem \eqref{bound_z1}, \eqref{bound_z2} is understood as a vector-valued function $y \in \big(C^{(n+r)}\big)^m$ that satisfies the differential equation \eqref{bound_z1} everywhere on $(a,b)$ as well as the boundary condition \eqref{bound_z2}.
As noted above, we refer to the general boundary condition \eqref{bound_z2} with an arbitrary continuous operator \eqref{oper_B(e)v} as \textit{generic} for the differential system \eqref{bound_z1}. It encompasses all classical types of boundary conditions (such as Cauchy initial conditions, multipoint conditions, integral conditions, and mixed conditions), as well as non-classical conditions containing integer and/or fractional derivatives of orders $s$, where $r\leq s \leq n+r$. Finally, we impose no {\it a priori} assumptions on the regularity of the matrix-valued functions $A_{r-\ell}(t,\mu)$ with respect to $\mu$.

We represent problem \eqref{bound_z1}, \eqref{bound_z2} in the form of an operator equation
\[ \big(L(\mu),B(\mu)\big)y(\mu)=\big(f(\mu),c(\mu)\big), \]
where $\big(L(\mu),B(\mu)\big)$ is the family of continuous linear operators
\begin{equation}\label{(L,B)vp}
\big(L(\mu),B(\mu)\big) \colon (C^{(n+r)})^m\to (C^{(n)})^m \oplus \mathbb{C}^{rm}.
\end{equation}

Let $E_{1}$ and $E_{2}$ be Banach spaces. Recall that a linear operator $T\colon E_{1}\rightarrow E_{2}$ is called Fredholm if its kernel $\ker \, T:= \{x \in E_{1} \colon Tx=0 \}$ and cokernel $\mathrm{coker}\, T:= E_{2}/ T (E_{1})$ are finite-dimensional. If $T$ is Fredholm, then its image $T(E_{1})$ is closed in $E_{2}$, and the integer
$$
\mathrm{ind}\,T:=\dim\ker T-\dim\big(E_{2}/T(E_{1})\big)\in \mathbb{Z} = \{0, \pm 1, \pm2 , \ldots \}
$$
is called its index (see, e.g., \cite[Sect.~19.1.1]{Hermander1985}).
According to \cite[Theorem 1]{Soldatov2025}, the family of operators \eqref{(L,B)vp} consists of Fredholm operators of index zero for each $\mu$. The operator \eqref{(L,B)vp} is a finite-dimensional perturbation of the operator acting from $(C^{(n+r)})^m$ to $(C^{(n)})^m\oplus \mathbb{C}^{rm}$ defined by the Cauchy problem for equation \eqref{bound_z1} (see, e.g., \cite{Hermander1985}).

\begin{definition}\label{defin_vp}
Problem \eqref{bound_z1}, \eqref{bound_z2} is said to be continuously dependent on the parameter $\mu$ at the point $\mu_0 \in \mathcal{M}$ if the following conditions are satisfied:
\begin{itemize}
\item[$(\ast)$] there exists a positive number $\varepsilon$ such that for each $\mu\in \mathcal{B}(\mu_0, \varepsilon)$ and any right-hand sides $f(\cdot;\mu)\in (C^{(n)})^{m}$ and $c(\mu)\in\mathbb{C}^{rm}$, the problem has a unique solution $y(\cdot;\mu) \in (C^{(n+r)})^{m}$;
\item [$(\ast\ast)$] the convergence of right-hand sides $f(\cdot;\mu)\to f(\cdot;\mu_0)$ in $(C^{(n)})^{m}$ and $c(\mu)\to c(\mu_0)$ in $\mathbb{C}^{rm}$ as $\mu\to\mu_0$ implies the convergence of solutions
\begin{equation*}\label{4.guv}
y(\cdot,\mu)\to y(\cdot,\mu_0)\quad\mbox{in}\quad (C^{(n+r)})^{m} \quad\mbox{as}\quad\mu\to\mu_0
\end{equation*}
for any $f(\cdot;\mu_0)\in (C^{(n)})^{m}$ and $c(\mu_0)\in\mathbb{C}^{rm}$.
\end{itemize}
\end{definition}

We now formulate the limit conditions at $\mu_0 \in \mathcal{M}$.

\medskip

\noindent \textbf{Condition (0)}. {\it The homogeneous limit boundary-value problem
\begin{equation*}
L(\mu_0)y(t,\mu_0)=0, \quad  t\in (a,b), \quad B(\mu_0)y(\cdot,\mu_0)=0
\end{equation*}
has only the trivial solution $y(t,\mu_0)\equiv 0$, $t\in (a,b)$ in $(C^{(n+r)})^m$. }

\medskip

\noindent We also introduce the following limit conditions for the problem \eqref{bound_z1}, \eqref{bound_z2}:

\medskip

\noindent \textbf{Limit Conditions as $\mu\to\mu_0$}:

\smallskip

\noindent
(I) {\it $A_{r-\ell}(\cdot;\mu)\to A_{r-\ell}(\cdot;\mu_0)$ in $(C^{(n)})^{m\times m}$ for each $\ell\in\{1,\ldots, r\}$;}

\smallskip

\noindent
(II) {\it $B(\mu)y\to B(\mu_0)y$ in $\mathbb{C}^{rm}$ for every $y\in(C^{(n+r)})^m$.
}

\medskip
The following theorem provides a criterion for the continuous parameter dependence of problem \eqref{bound_z1}, \eqref{bound_z2}.

\begin{theorem}\label{nep v}
Problem \eqref{bound_z1}, \eqref{bound_z2} is continuously dependent on the parameter $\mu$ at $\mu_0\in \mathcal{M}$ if and only if Condition~\textup{(0)} and the Limit Conditions~\textup{(I)} and~\textup{(II)} hold.
\end{theorem}

\begin{corollary}
If Condition~\textup{(0)} and Limit Conditions~\textup{(I)} and~\textup{(II)} are satisfied for every $\mu \in \mathcal{M}$, then for any continuous vector-valued functions $f$ and $c$, the solution of problem \eqref{bound_z1}, \eqref{bound_z2} depends continuously on $\mu \in \mathcal{M}$ as an element of the space $C\big(\mathcal{M}; (C^{(n+r)})^{m}\big)$.
\end{corollary}

Notice that the metric space $\mathcal{M}$ in Theorem~\ref{nep v} is arbitrary; it need not be locally compact or separable.

Prior to proving Theorem~\ref{nep v}, we establish an auxiliary result, which is an analogue of \cite[Theorem~3.4]{MikhMurachSol2016} for the case $\alpha=0$ in the spaces under consideration.

We introduce the metric space
$$ \mathcal{Y}^{n+1}:=\{Y(\cdot)\in (C^{(n+1)})^{m\times m}\colon \,Y(a)=I_{m},\quad
  \det Y(t)\neq 0, \quad t\in [a, b]\}$$
\noindent equipped with the metric
 $$d_{n+1}(Y,Z):=\|Y(\cdot)-Z(\cdot)\|_{(n+1)}.$$

\begin{theorem}\label{th1}
The mapping
\begin{equation*}\label{rr6}
A(\cdot)\mapsto Y(\cdot),
\end{equation*}
which associates with each matrix-valued function $A(\cdot) \in (C^{(n)})^{m\times m}$ the solution $Y(\cdot)$ of the Cauchy problem
\begin{equation*}\label{r3}
Y'(t)+A(t)Y(t)=0,\quad t\in (a,b), \quad Y(a)=I_{m},
\end{equation*}
is a homeomorphism of the space $(C^{(n)})^{m\times m}$ onto $\mathcal{Y}^{n+1}$.
\end{theorem}

\begin{proof}[Proof of Theorem \ref{nep v}]

First, we prove the \emph{necessity} part of Theorem~\ref{nep v}. Suppose that problem \eqref{bound_z1}, \eqref{bound_z2} is continuously dependent on $\mu$ in the sense of Definition~\ref{defin_vp}.
Condition~(0) is clearly satisfied by virtue of condition~$(\ast)$. We must demonstrate the necessity of conditions~(I) and~(II). We divide the argument into three steps.

\emph{Step 1.} We prove Condition~\textup{(I)}. It is well known that an $r$-th order differential system can be reduced to a first-order system of dimension $rm$ (see, e.g., \cite{Cartan1971}). Denoting the transpose by $\top$, we set
\begin{gather}\label{yx}
x(\cdot,\mu):= \bigl(y(\cdot,\mu),y'(\cdot, \mu),\ldots,y^{(r-1)}(\cdot, \mu)\bigr)^\top\in(C^{(n+1)})^{rm},\\
g(\cdot,\mu):= \bigl(\underbrace{0, \dots, 0}_{(r-1)m},f(\cdot, \mu)\bigr)^\top
\in(C^{(n)})^{rm},\notag\\
c(\mu):= \bigl(c_{1}(\mu),\ldots,c_{r}(\mu)\bigr)^\top
\in\mathbb{C}^{rm}\notag,
\end{gather}
and define the matrix-valued function $K(\cdot,\mu) \in (C^{(n)})^{rm\times rm}$ by
\begin{equation}\label{AAv-reduced}
K(\cdot,\mu):=\left(
\begin{array}{ccccc}
O_m & -I_m & O_m & \ldots & O_m \\
O_m & O_m & -I_m & \ldots & O_m \\
\vdots & \vdots & \vdots & \ddots & \vdots \\
O_m & O_m & O_m & \ldots & -I_m \\
A_0(\cdot,\mu) & A_1(\cdot,\mu) & A_2(\cdot,\mu) & \ldots & A_{r-1}(\cdot,\mu)\\
\end{array}\right),
\end{equation}
where $O_m$ and $I_m$ denote the zero and identity matrices of dimension $m \times m$, respectively.
A vector-valued function $y(\cdot,\mu)\in(C^{(n+r)})^{m}$ is a solution of system \eqref{bound_z1} if and only if the vector-valued function \eqref{yx} is a solution of the system
\begin{gather*}
x'(t,\mu)+K(t,\mu)x(t,\mu)=g(t,\mu),
\quad t\in(a,b).
\end{gather*}

We denote by
\begin{equation*}\label{matrix_BYe}
\left[B(\mu)Y(\cdot,\mu)\right]:=\left(\left[B(\mu)Y_0(\cdot,\mu)\right],\dots,\left[B(\mu)
Y_{r-1}(\cdot,\mu)\right]\right) \in \mathbb{C}^{rm\times rm}
\end{equation*}
the block matrix of dimension $rm\times rm$. It consists of $r$ block columns $\left[B(\mu)Y_k(\cdot,\mu)\right]\in \mathbb{C}^{rm\times m}$ ($k \in \{0, \dots, r-1\}$), where for each $s \in \{1, \dots, m\}$ the $s$-th column of the matrix $\left[B(\mu)Y_k(\cdot,\mu)\right]$ is the result of the action of the operator $B(\mu)$ on the $s$-th column of the matrix-valued function $Y_k(\cdot,\mu)$.

Consider the matrix-valued function $Y(\cdot, \mu) := (Y_0(\cdot, \mu), \ldots, Y_{r-1}(\cdot, \mu)) \in (C^{(n+r)})^{m \times rm}$, composed of $r$ blocks of dimension $m \times m$:
\begin{equation}\label{eq2vk}
\begin{gathered}
  Y^{(r)}(t, \mu) + \sum_{\ell=1}^r A_{r-\ell}(t, \mu) Y^{(r-\ell)}(t, \mu) = O_{m \times rm}, \quad t \in (a, b), \\
  [B(\mu)Y(\cdot, \mu)] = I_{rm}.
\end{gathered}
\end{equation}
Here $O_{m \times rm}$ and $I_{rm}$ denote the zero matrix of dimension $m \times rm$ and the identity matrix of dimension $rm \times rm$, respectively. System \eqref{eq2vk} represents $rm$ problems of the form \eqref{bound_z1}, \eqref{bound_z2} with zero right-hand sides for $f$, which are independent of $\mu$.

Thus, there exists a unique solution $Y(\cdot;\mu)$ for each $\mu\in \mathcal{B}(\mu_0, \varepsilon)$, according to property~$(\ast)$ of Definition~\ref{defin_vp}. Furthermore, condition~$(\ast\ast)$ implies
\begin{equation}\label{zb ym}
y_{j,k}(\cdot,\mu)\rightarrow y_{j,k}(\cdot,\mu_0) \quad\mbox{in} \quad C^{(n+r)} \quad\mbox{as}\quad \mu\to\mu_0.
\end{equation}

For each $k\in\{1, \ldots, rm\}$ and $\mu\in \mathcal{B}(\mu_0, \varepsilon)$, we define the vector-valued function $x_k(\cdot,\mu)\in(C^{(n+1)})^{rm}$ by formula \eqref{yx}, where we set
$$
y(\cdot;\mu):=\big(y_{1,k}(\cdot;\mu),\ldots, y_{m,k}(\cdot;\mu)\big)^\top.
$$
Let $X(\cdot;\mu) \in (C^{(n+1)})^{rm\times rm}$ denote the matrix-valued function whose $k$-th column is $x_k(\cdot,\mu)$ for each $k\in\{1, \ldots, rm\}$. This function satisfies the matrix differential equation
\begin{equation}\label{X ym}
X'(t,\mu)+K(t,\mu)X(t,\mu)=O_{rm}, \quad t\in (a,b),
\end{equation}
where $O_{rm}$ is the zero matrix of dimension $rm \times rm$.
Since $\det X(t;\mu)\neq0$ for all $t\in[a,b]$, the columns of $X(\cdot;\mu)$, and consequently the columns of $Y(\cdot;\mu)$, form a fundamental system of solutions on $[a,b]$, satisfying \eqref{eq2vk}.
Relation \eqref{zb ym} yields $X(\cdot;\mu)\rightarrow X(\cdot;\mu_0)$ in $(C^{(n+1)})^{rm\times rm}$ as $\mu\to\mu_0$. Hence, $\big(X(\cdot;\mu)\big)^{-1}\rightarrow \big(X(\cdot;\mu_0)\big)^{-1}$ in the same space. Using equation \eqref{X ym}, we obtain
  $$
   K(\cdot;\mu)=-X'(\cdot;\mu)\big(X(\cdot;\mu)\big)^{-1} \rightarrow -X'(\cdot;\mu_0)\big(X(\cdot;\mu_0)\big)^{-1}=K(\cdot;\mu_0)
  $$
in $(C^{(n)})^{rm\times rm}$ as $\mu\to\mu_0$. In view of the structure of $K$ in \eqref{AAv-reduced}, this implies that problem \eqref{bound_z1}, \eqref{bound_z2} satisfies Limit Condition~(I). In particular,
\begin{equation}\label{ob Av}
\big\|A_{r-\ell}(\cdot,\mu)\big\|_{(n)}=O(1) \quad \mbox{as}\ \mu\to\mu_0,  \ \mbox{for each} \ \ell\in\{1, \ldots, r\}.
\end{equation}

\emph{Step 2.} We now show that
\begin{equation}\label{ob Bv}
\|B(\mu)\|=O(1) \quad\mbox{as}\ \mu\to\mu_0.
\end{equation}

Suppose, on the contrary, that this is false; then there exists a sequence $\left(\mu^{(k)}\right)_{k=1}^{\infty}\subset\mathcal{B}(\mu_0, \varepsilon)$ such that
\begin{equation}\label{ob Bev}
\mu^{(k)}\to\mu_0\quad\mbox{and}\quad 0<\left\|B\bigl(\mu^{(k)}\bigr)\right\|\to\infty \quad\mbox{as}
\ k\rightarrow \infty,
\end{equation}
where $\left\|B\bigl(\mu^{(k)}\bigr)\right\|\neq0$ for all $k\in \mathbb{N}$. For each $k\in \mathbb{N}$, we choose a vector-valued function $\omega_{k}\in(C^{(n+r)})^{m}$ such that
\begin{equation}\label{ob wv}
\|\omega_{k}\|_{(n+r)}=1 \quad \mbox{and} \quad \bigl\|B\bigl(\mu^{(k)}\bigr)\omega_{k}\bigr\|_{\mathbb{C}^{rm}}\geq \frac{1}{2}
\bigl\|B\bigl(\mu^{(k)}\bigr)\bigr\|.
\end{equation}
We define
\begin{gather*}
y\bigl(\cdot;\mu^{(k)}\bigr):=
\bigl\|B\bigl(\mu^{(k)}\bigr)\bigr\|^{-1}\omega_{k} \in (C^{(n+r)})^{m},\\
f\bigl(\cdot;\mu^{(k)}\bigr):=
L\bigl(\mu^{(k)}\bigr)\,y\bigl(\cdot;\mu^{(k)}\bigr) \in (C^{(n)})^{m},\\
c\bigl(\mu^{(k)}\bigr):=B\bigl(\mu^{(k)}\bigr)\,y\bigl(\cdot;\mu^{(k)}\bigr) \in \mathbb{C}^{rm}. \end{gather*}
From \eqref{ob Bev} and \eqref{ob wv}, it follows that
\begin{equation}\label{zb ymal}
y\left(\cdot;\mu^{(k)}\right)\to0 \quad \mbox{in} \ (C^{(n+r)})^{m}  \  \mbox{as} \ k\rightarrow \infty.
\end{equation}
Consequently,
\begin{equation}\label{zb fmal}
f\left(\cdot;\mu^{(k)}\right)\to0 \quad \mbox{in}\ (C^{(n)})^{m}  \  \mbox{as} \ k\rightarrow \infty,
\end{equation}
since problem \eqref{bound_z1}, \eqref{bound_z2} satisfies Condition~(I) by Step~1.

Moreover, because the space $\mathbb{C}^{rm}$ is finite-dimensional, relation \eqref{ob wv} implies
$$
1/2\leq\left\|c\left(\mu^{(k)}\right)\right\|_{\mathbb{C}^{rm}}\leq1.
$$
Indeed,
\begin{gather*}
\left\|c\left(\mu^{(k)}\right)\right\|_{\mathbb{C}^{rm}}\leq \left\|B\left(\mu^{(k)}\right)\right\| \left\|y\left(\cdot,\mu^{(k)}\right)\right\|_{(n+r)}=\\
=\left\|B\left(\mu^{(k)}\right)\right\| \left\|B\left(\mu^{(k)}\right)\right\|^{-1} \left\|\omega_k\right\|_{(n+r)}=1
\end{gather*}
and
\begin{gather*}
\left\|c\left(\mu^{(k)}\right)\right\|_{\mathbb{C}^{rm}}= \left\|B\left(\mu^{(k)}\right)\left(
\left\|B\left(\mu^{(k)}\right)\right\|^{-1}\omega_k\right)\right\| = \left\|B\left(\mu^{(k)}\right)\right\|^{-1} \left\|B\left(\mu^{(k)}\right)\omega_k\right\|\geq \frac{1}{2}.
\end{gather*}
Thus, there exists a subsequence $$\left(c\left(\mu^{(k_s)}\right)\right)^\infty_{s=1}\subset \left(c\left(\mu^{(k)}\right)\right)^\infty_{k=1}$$ and a non-zero vector $c(\mu_0)\in \mathbb{C}^{rm}$ such that
\begin{equation}\label{zb cmal}
c\left(\mu^{(k_s)}\right)\to c(\mu_0) \quad \mbox{in} \ \mathbb{C}^{rm}  \  \mbox{as} \ s\rightarrow \infty.
\end{equation}
For each $s\in \mathbb{N}$, the function $y\left(\cdot;\mu^{(k_s)}\right)\in\left(C^{(n+r)}\right)^{m}$ is the unique solution of the problem
\begin{gather*}
L\left(\mu^{(k_s)}\right)y\left(t;\mu^{(k_s)}\right)=f\left(t;\mu^{(k_s)}\right),  \quad t\in (a,b),\\
B\left(\mu^{(k_s)}\right)y\left(\cdot;\mu^{(k_s)}\right)=c\left(\mu^{(k_s)}\right).
\end{gather*}
By virtue of \eqref{zb fmal} and \eqref{zb cmal}, together with condition~$(\ast\ast)$ of Definition~\ref{defin_vp}, we conclude that $y\left(\cdot;\mu^{(k_s)}\right)$ converges to the unique solution $y(\cdot;\mu_0)$ of the limit problem
\begin{gather}
L(\mu_0)y(t,\mu_0)=0,\quad t\in (a,b), \nonumber\\
B(\mu_0)y(\cdot;\mu_0)=c(\mu_0) \label{kr ym}
\end{gather}
in the space $(C^{(n+r)})^{m}$ as $s\rightarrow \infty$. However, $y(\cdot;\mu_0)\equiv0$ due to \eqref{zb ymal}. This contradicts the boundary condition \eqref{kr ym}, in which $c(\mu_0)\neq0$. This proves assertion \eqref{ob Bv}.

\emph{Step 3.} Using the results of the preceding steps, we now show that problem \eqref{bound_z1}, \eqref{bound_z2} satisfies Limit Condition~(II). According to \eqref{ob Av} and \eqref{ob Bv}, there exist constants $\gamma'>0$ and $\varepsilon'>0$ such that
\begin{equation}\label{ner 1}
\|(L(\mu),B(\mu))\|\leq\gamma' \quad \mbox{for all}  \  \mu\in\mathcal{B}(\mu_0, \varepsilon'),
\end{equation}
where $\|\cdot\|$ denotes the operator norm on \eqref{(L,B)vp}. For an arbitrary vector-valued function $y\in(C^{(n+r)})^{m}$, we set $f(\cdot;\mu):=L(\mu)y$ and $c(\mu):=B(\mu)y$. Then
\begin{equation}\label{riv 1}
y=(L(\mu),B(\mu))^{-1}(f(\cdot;\mu),c(\mu)) \quad \mbox{for all}  \  \mu\in\mathcal{B}(\mu_0, \varepsilon').
\end{equation}
Here $(L(\mu),B(\mu))^{-1}$ is the inverse operator to \eqref{(L,B)vp}, which exists by property~$(\ast)$ of Definition~\ref{defin_vp}.
From \eqref{ner 1} and \eqref{riv 1}, we introduce the difference
\[
w(\mu) := (f(\cdot;\mu), c(\mu)) - (f(\cdot;\mu_0), c(\mu_0)) \in (C^{(n)})^m \oplus \mathbb{C}^{rm}.
\]
Since $w(\mu) = (L(\mu), B(\mu))\,(L(\mu), B(\mu))^{-1} w(\mu)$, for all $\mu\in\mathcal{B}(\mu_0, \varepsilon')$ we have:
\begin{align*}
\bigl\|B(\mu)y - B(\mu_0)y\bigr\|_{\mathbb{C}^{rm}}
&\leq \|w(\mu)\|_{(C^{(n)})^m \times \mathbb{C}^{rm}} \\
&= \bigl\| (L(\mu), B(\mu))\,(L(\mu), B(\mu))^{-1} w(\mu) \bigr\|_{(C^{(n)})^m \times \mathbb{C}^{rm}} \\
&\leq \gamma'\, \bigl\| (L(\mu), B(\mu))^{-1} w(\mu) \bigr\|_{(n+r)}.
\end{align*}
Furthermore, using the identity \eqref{riv 1}:
\[
(L(\mu), B(\mu))^{-1}(f(\cdot;\mu), c(\mu)) = y = (L(\mu_0), B(\mu_0))^{-1}(f(\cdot;\mu_0), c(\mu_0)),
\]
we observe that
\[
(L(\mu), B(\mu))^{-1} w(\mu) = (L(\mu_0), B(\mu_0))^{-1}(f(\cdot;\mu_0), c(\mu_0)) - (L(\mu), B(\mu))^{-1}(f(\cdot;\mu_0), c(\mu_0)).
\]
This expression converges to $0$ in $(C^{(n+r)})^m$ as $\mu\to\mu_0$ by property~$(\ast\ast)$ of Definition~\ref{defin_vp}.

Therefore, relation \eqref{ob Bv} implies that problem \eqref{bound_z1}, \eqref{bound_z2} satisfies Limit Condition~(II), which completes the proof of necessity.

\medskip

We now turn to the proof of \emph{sufficiency} in Theorem~\ref{nep v}. Suppose that problem \eqref{bound_z1}, \eqref{bound_z2} satisfies Condition~(0) as well as Limit Conditions~(I) and~(II). We must demonstrate that the solution of the problem exists, is unique in $(C^{(n+r)})^{m}$ for all $\mu$ sufficiently close to $\mu_0$, and satisfies the continuity condition~$(\ast\ast)$. We divide the argument into four steps.

\emph{Step 1.} For each $\mu\in \mathcal{M}$, we first consider the Cauchy initial-value problem with initial conditions at an arbitrary fixed point $t_0 \in [a, b]$ (without loss of generality, we set $t_0 = a$):
\begin{gather}\label{s1}
L(\mu)\hat{y}(t,\mu)=f(t,\mu),\quad t\in(a,b),\\
\hat{y}^{(j-1)}(a,\mu)=c_{j}(\mu),\quad j\in\{1,\ldots,r\}. \label{ku1}
\end{gather}
Here, for each $\mu$, the function $f(\cdot,\mu)\in(C^{(n)})^m$ and the vectors $c_{j}(\mu)\in\mathbb{C}^{m}$ are given arbitrarily. The unique solution $\hat{y}(\cdot,\mu)$ of this Cauchy problem belongs to $(C^{(n+r)})^{m}$.

We prove that the convergence of the right-hand sides
\begin{equation}\label{r_s.1}
f(\cdot,\mu)\to f(\cdot,\mu_0)\quad\mbox{in} \ (C^{(n)})^{m}  \ \mbox{as} \  \mu\to \mu_0,
\end{equation}
\begin{equation}\label{gsp}
\begin{gathered}
c_{j}(\mu)\rightarrow c_{j}(\mu_0)\quad \mbox{in}  \  \mathbb{C}^{m} \ \mbox{as} \  \mu\to \mu_0
 \ \mbox{for each} \  j\in\{1,\ldots,r\}
\end{gathered}
\end{equation}
implies the convergence of the solutions
\begin{equation}\label{gsx}
\hat{y}(\cdot,\mu)\to \hat{y}(\cdot,\mu_0)\quad\mbox{in} \ (C^{(n+r)})^{m}  \ \mbox{as} \  \mu\rightarrow\mu_0.
\end{equation}
The Cauchy problem \eqref{s1}, \eqref{ku1} can be reduced to the following initial-value problem for a first-order system:
\begin{gather}\label{Cauchi-reduced-DE}
x'(t,\mu)+K(t,\mu)x(t,\mu)=g(t,\mu),
\quad t\in(a,b),\\
x(a,\mu)=c(\mu).\label{Cauchi-reduced-cond}
\end{gather}
Here the matrix-valued function $K(\cdot,\mu)$ and vector-valued function $g(\cdot,\mu)$ are defined as in Step~1 of the necessity proof. In particular,
\begin{gather*}\label{ghggh}
x(\cdot,\mu):=\left(\hat{y}(\cdot,\mu),\hat{y}'(\cdot,\mu),\ldots,\hat{y}^{(r-1)}(\cdot,\mu)\right)^\top\in (C^{(n+1)})^{rm}, \\
c(\mu):= \left(c_1(\mu),\ldots,c_r(\mu)\right)^\top\in \mathbb{C}^{rm}.
\end{gather*}
Since problem \eqref{bound_z1}, \eqref{bound_z2} satisfies Limit Condition~(I), we have
\begin{equation*}\label{Kgsx}
K(\cdot,\mu)\to K(\cdot,\mu_0)\quad\mbox{in} \ (C^{(n)})^{rm\times rm}  \ \mbox{as} \  \mu\rightarrow\mu_0.
\end{equation*}
Moreover, relations \eqref{r_s.1} and \eqref{gsp} ensure that for the problem \eqref{Cauchi-reduced-DE}, \eqref{Cauchi-reduced-cond} we have
\begin{equation*}\label{r_g.1}
g(\cdot,\mu)\to g(\cdot,\mu_0)\quad\mbox{in} \ (C^{(n)})^{rm}  \ \mbox{as} \  \mu\to\mu_0,
\end{equation*}
\begin{equation*}\label{ggp}
\begin{gathered}
c(\mu)\rightarrow c(\mu_0)\quad\mbox{in} \ \mathbb{C}^{rm} \ \mbox{as} \  \mu\to \mu_0.
\end{gathered}
\end{equation*}
Hence, by \cite[Theorem~3]{Sold2015}, relations \eqref{r_s.1} and \eqref{gsp} indeed imply \eqref{gsx}.

\emph{Step 2.} We now prove that problem \eqref{bound_z1}, \eqref{bound_z2} satisfies condition~$(\ast)$ of Definition~\ref{defin_vp}, that is, for sufficiently small $d(\mu,\mu_0)$ the operator $(L(\mu),B(\mu))$ is invertible.

For each $\mu\in \mathcal{M}$ and $k\in\{0,\ldots,r-1\}$, consider the matrix Cauchy problem:
\begin{gather}\label{zad kosh1v}
Y_k^{(r)}(t,\mu)+\sum\limits_{j=1}^rA_{r-j}(t,\mu)Y_k^{(r-j)}(t,\mu)=O_{m},\quad t\in (a,b),
\end{gather}
with the initial conditions
\begin{gather*}\label{zad kosh2v}
Y_k^{(j)}(t_0,\mu) = \delta_{k,j}I_m,\quad j \in \{0,\dots, r-1\}.
\end{gather*}
Here
$$
Y_k(t,\mu)=\left(y_k^{\alpha,\beta}(t,\mu)\right)_{\alpha,\beta=1}^m
$$
is an unknown $m\times m$ matrix-valued function, the point $t_0 \in [a,b]$ is fixed (without loss of generality, $t_0 = a$), $\delta_{k,j}$ is the Kronecker delta, and $O_{m}$ and $I_m$ are as in \eqref{AAv-reduced}. Problem \eqref{zad kosh1v} consists of $m$ Cauchy problems of the form \eqref{s1}, \eqref{ku1} with $f=0$ for the unknown vector-valued functions $\hat{y}(\cdot,\mu)$, which are the columns of $Y_k(\cdot,\mu)$.

We express the general solution of the homogeneous differential equation \eqref{bound_z1} in the form
\begin{equation}\label{solution_by_Yv}
y(\cdot,\mu)=\sum\limits_{k=0}^{r-1}Y_k(\cdot,\mu)q_k(\mu),
\end{equation}
where $q_k(\mu)\in\mathbb{C}^{m}$ are arbitrary column vectors \cite{Cartan1971}.
The right-hand sides of this problem do not depend on $\mu$; therefore,
\begin{equation}\label{Zlimitv}
Y_k(\cdot,\mu)\to Y_k(\cdot,\mu_0) \quad\mbox{in} \ (C^{(n+r)})^{m\times m} \ \mbox{as} \  \mu\to\mu_0,
\end{equation}
according to Step~1 and Theorem~\ref{th1}.
To justify the convergence of the matrices $[B(\mu)Y_k(\cdot,\mu)]$, we use the identity
\begin{equation*}
B(\mu)Y_k(\cdot,\mu) - B(\mu_0)Y_k(\cdot,\mu_0) = B(\mu)\bigl(Y_k(\cdot,\mu) - Y_k(\cdot,\mu_0)\bigr) + \bigl(B(\mu) - B(\mu_0)\bigr)Y_k(\cdot,\mu_0).
\end{equation*}
By Limit Condition~(II), the operators $B(\mu)$ converge strongly to $B(\mu_0)$ as $\mu\to\mu_0$. By the Banach--Steinhaus theorem, their norms are uniformly bounded: there exists a constant $C > 0$ such that $\|B(\mu)\| \le C$ for all $\mu$ in a neighborhood of $\mu_0$. Consequently,
\begin{equation*}
\|B(\mu)\bigl(Y_k(\cdot,\mu) - Y_k(\cdot,\mu_0)\bigr)\| \le C \|Y_k(\cdot,\mu) - Y_k(\cdot,\mu_0)\|_{(n+r)} \to 0,
\end{equation*}
while the second term $\bigl(B(\mu) - B(\mu_0)\bigr)Y_k(\cdot,\mu_0) \to 0$ directly by Condition~(II). This yields the convergence
\begin{equation}\label{BZlimit}
\begin{gathered}
\bigl([B(\mu)Y_0(\cdot,\mu)],\ldots,
[B(\mu)Y_{r-1}(\cdot,\mu)]\bigr)\to\\
\to\bigl([B(\mu_0)Y_0(\cdot,\mu_0)],\ldots,
[B(\mu_0)Y_{r-1}(\cdot,\mu_0)]\bigr)\quad \mbox{in}  \  \mathbb{C}^{rm\times rm} \
\mbox{as} \ \mu\to\mu_0.
\end{gathered}
\end{equation}
Furthermore, the limit matrix is non-singular by virtue of Condition~(0) and \cite[Theorems~1, 2]{Sold2015}. Hence, there exists $\varepsilon_1 >0$ such that for each $\mu\in\mathcal{B}(\mu_0, \varepsilon_1)$,
\begin{equation}\label{BZdet}
\begin{gathered}
\det\bigl(M(L(\mu),B(\mu))\bigr)\neq0.
\end{gathered}
\end{equation}
Thus, according to \cite[Theorem~2]{Sold2015}, the operator \eqref{(L,B)vp} is invertible for every $\mu\in\mathcal{B}(\mu_0, \varepsilon_1)$.

\emph{Step 3.} We now show that problem \eqref{bound_z1}, \eqref{bound_z2} satisfies condition~$(\ast\ast)$ of Definition~\ref{defin_vp}. First, we analyze the case $f(\cdot,\mu)\equiv0$: consider the $\mu$-dependent semi-homogeneous boundary-value problem
 \begin{equation}\label{vv}
    L(\mu)v(\cdot;\mu)\equiv 0,
\end{equation}
 \begin{equation}\label{vvk}
  B(\mu)v(\cdot;\mu) =c(\mu).
\end{equation}

By Step~2, this problem has a unique solution $v(\cdot;\mu)\in(C^{(n+r)})^{m}$ for each $\mu\in \mathcal{B}(\mu_0, \varepsilon_1)$. Assuming that
 \begin{equation}\label{zb c2v}
  c(\mu)\rightarrow c(\mu_0) \quad \mbox{in}  \  \mathbb{C}^{rm}  \  \mbox{as}  \ \mu\rightarrow \mu_0,
 \end{equation}
we show that
 \begin{equation}\label{zb v2v}
  v(\cdot,\mu)\rightarrow v(\cdot,\mu_0) \quad \mbox{in}  \  (C^{(n+r)})^{m}  \  \mbox{as}  \ \mu\rightarrow \mu_0.
 \end{equation}

For $\mu\in\mathcal{B}(\mu_0, \varepsilon_1)$, the solution of \eqref{vv} can be represented in the form \eqref{solution_by_Yv}:
\begin{equation}\label{solution_by_vv}
v(\cdot,\mu)=\sum\limits_{k=0}^{r-1}Y_k(\cdot,\mu)q_k(\mu),
\end{equation}
for some vectors $q_k(\mu)\in\mathbb{C}^{m}$ and matrix-valued functions $Y_k(\cdot,\mu)\in(C^{(n+r)})^{m\times m}$, as in Step~2. By \cite[Lemma~5.1]{NovAM2023},
\begin{equation*}
B(\mu)v(\cdot,\mu)=\sum _{k=0}^{r-1}B(\mu)(Y_k(\cdot,\mu)q_k(\mu)) = \sum _{k=0}^{r-1}\left[B(\mu)Y_k(\cdot,\mu)\right]q_k(\mu),
\end{equation*}
so the boundary condition \eqref{vvk} is equivalent to
\begin{equation*}
\sum _{k=0}^{r-1}\left[B(\mu)Y_k(\cdot,\mu)\right]q_k(\mu)=c(\mu).
\end{equation*}
This equation can be rewritten as the linear algebraic system
\begin{equation*}
\big(\left[B(\mu)Y_0(\cdot,\mu)\right], \ldots, \left[B(\mu)Y_{r-1}(\cdot,\mu)\right]\big)q(\mu)=c(\mu),
 \ \mbox{or}  \
\left[B(\mu)Y(\cdot,\mu)\right]q(\mu)=c(\mu) ,
\end{equation*}
where $q(\mu):=(q_0(\mu),\ldots, q_{r-1}(\mu))^\top.$
Now, relations \eqref{BZlimit}, \eqref{BZdet}, and assumption \eqref{zb c2v} yield
\begin{equation*}
q(\mu)=\left[B(\mu)Y(\cdot,\mu)\right]^{-1}c(\mu)\rightarrow\left[B(\mu_0)Y(\cdot,\mu_0)\right]^{-1}c(\mu_0)=q(\mu_0) \quad \mbox{as}  \ \mu\to \mu_0.
\end{equation*}
Consequently, \eqref{zb v2v} follows from \eqref{Zlimitv} and \eqref{solution_by_vv}, since
\begin{equation*}
v(\cdot,\mu)=\sum _{k=0}^{r-1}Y_k(\cdot,\mu)q_k(\mu) \rightarrow \sum _{k=0}^{r-1}Y_k(\cdot,\mu_0)q_k(\mu_0)=v(\cdot,\mu_0)
\end{equation*}
in $(C^{(n+r)})^{m}$ as $\mu\to \mu_0.$

\emph{Step 4.} Finally, we consider the general case of the inhomogeneous differential equation \eqref{bound_z1}. Suppose that condition \eqref{zb c2v} and
\begin{equation}\label{zb fvvy}
f(\cdot,\mu) \rightarrow f(\cdot,\mu_0) \quad \mbox{in}  \  (C^{(n)})^{m}  \  \mbox{as}  \ \mu\to \mu_0
\end{equation}
are fulfilled.
For each $\mu \in \mathcal{B}(\mu_0, \varepsilon_1)$, we set
$$z(\cdot;\mu)=y(\cdot;\mu)-\hat{y}(\cdot;\mu),$$
where the vector-valued function $y(\cdot;\mu)$ is the solution of the inhomogeneous boundary-value problem \eqref{bound_z1}, \eqref{bound_z2}, and $\hat{y}(\cdot;\mu)$ is the solution of the Cauchy problem \eqref{s1}, \eqref{ku1} with $c_{j}(\mu)\equiv0$. Then $z(\cdot;\mu)$ is the solution of the semi-homogeneous boundary-value problem
\begin{gather*}\label{s2}
L(\mu)z(\cdot,\mu)\equiv0,\\
B(\mu)z(\cdot,\mu)=\tilde{c}(\mu),\\
\tilde{c}(\mu)=c(\mu)-B(\mu)\hat{y}(\cdot,\mu) \in \mathbb{C}^{rm}.
\end{gather*}
In Step~1, it was established that $\hat{y}(\cdot,\mu)$ satisfies \eqref{gsx} whenever condition \eqref{zb fvvy} holds. This fact, together with \eqref{zb c2v}, the uniform boundedness $\|B(\mu)\| = O(1)$ from Step~2, and Limit Condition~(II), ensures that $\tilde{c}(\mu) \to \tilde{c}(\mu_0)$ in $\mathbb{C}^{rm}$ as $\mu \to \mu_0$.
Applying Step~3, we conclude that
\begin{equation*}
z(\cdot,\mu) \rightarrow z(\cdot,\mu_0) \quad \mbox{in}  \  (C^{(n+r)})^{m}  \  \mbox{as}  \ \mu\to \mu_0 ,
\end{equation*}
which yields the desired convergence from formula \eqref{gsx}:
\begin{gather*}
y(\cdot,\mu)= \hat{y}(\cdot,\mu)+z(\cdot,\mu)\rightarrow \hat{y}(\cdot,\mu_0)+z(\cdot,\mu_0)=y(\cdot,\mu_0)\\ \quad \mbox{in}  \  (C^{(n+r)})^{m}  \  \mbox{as}  \ \mu\to \mu_0.
\end{gather*}
This completes the proof of sufficiency.
\end{proof}

We now consider an approximation problem for the solution $y(\cdot;\mu)$ of problem \eqref{bound_z1}, \eqref{bound_z2} in terms of the unperturbed solution $y(\cdot;\mu_0)$. Specifically, we estimate the deviation
$$
\bigl\|y(\cdot;\mu_0)-y(\cdot;\mu)\bigr\|_{(n+r)}
$$
via the residual
\begin{equation*}\label{nevyuzka v}
\widetilde{d}_{(n)}(\mu):=
\bigl\|L(\mu)y(\cdot;\mu_0)-f(\cdot;\mu)\bigr\|_{(n)}+
\bigl\|B(\mu)y(\cdot;\mu_0)-c(\mu)\bigr\|_{\mathbb{C}^{rm}}.
\end{equation*}

\begin{theorem}\label{3.6.th-bound v}
Let $\mu_0 \in \mathcal{M}$ and suppose that problem \eqref{bound_z1}, \eqref{bound_z2} satisfies Condition~\textup{(0)} and Limit Conditions~\textup{(I)} and~\textup{(II)}. Then there exist positive numbers $\varepsilon$, $\gamma_{1}$, and $\gamma_{2}$ such that
\begin{equation}\label{3.6.bound}
\begin{aligned}
\gamma_{1}\,\widetilde{d}_{(n)}(\mu)
\leq\bigl\|y(\cdot;\mu_0)-y(\cdot;\mu)\bigr\|_{(n+r)}\leq
\gamma_{2}\,\widetilde{d}_{(n)}(\mu)
\end{aligned}
\end{equation}
for all $\mu\in\mathcal{B}(\mu_0, \varepsilon)$. Here the numbers $\varepsilon$, $\gamma_{1}$, and $\gamma_{2}$ do not depend on $y(\cdot;\mu_0)$ or $y(\cdot;\mu)$.
\end{theorem}

\begin{proof}
The left-hand inequality in \eqref{3.6.bound} follows directly from \eqref{ner 1}, choosing $\varepsilon$ no greater than $\varepsilon'$:
\begin{gather*}
\widetilde{d}_{(n)}(\mu)=\bigl\|L(\mu)y(\cdot;\mu_0)-L(\mu)y(\cdot;\mu)\bigr\|_{(n)}+
\bigl\|B(\mu)y(\cdot;\mu_0)-B(\mu)y(\cdot;\mu)\bigr\|_{\mathbb{C}^{rm}}\leq
\\
\leq\bigl\|L(\mu)\bigr\|\bigl\|y(\cdot;\mu_0)-y(\cdot;\mu)\bigr\|_{(n+r)}+
\bigl\|B(\mu)\bigr\|\bigl\|y(\cdot;\mu_0)-y(\cdot;\mu)\bigr\|_{(n+r)}
\leq
\\ \leq \gamma'\bigl\|y(\cdot;\mu_0)-y(\cdot;\mu)\bigr\|_{(n+r)}.
\end{gather*}

We now establish the right-hand inequality. By Theorem~\ref{nep v}, problem \eqref{bound_z1}, \eqref{bound_z2} satisfies the conditions of Definition~\ref{defin_vp}, so the operator \eqref{(L,B)vp} is invertible for all sufficiently small $d(\mu,\mu_0)$. In addition, the family of inverse operators satisfies
$$
(L(\mu),B(\mu))^{-1} \stackrel{s}{\longrightarrow} (L(\mu_0),B(\mu_0))^{-1}, \quad \mu\to\mu_0,
$$
where $\stackrel{s}{\longrightarrow}$ denotes strong operator convergence. Indeed, for any $f\in(C^{(n)})^{m}$ and $c\in\mathbb{C}^{rm}$, condition~$(\ast\ast)$ of Definition~\ref{defin_vp} gives
\begin{equation*}
(L(\mu),B(\mu))^{-1}(f,c)=:y(\cdot,\mu)\to
y(\cdot,\mu_0):=(L(\mu_0),B(\mu_0))^{-1}(f,c)
\end{equation*}
in $(C^{(n+r)})^{m}$ as $\mu\to\mu_0$.
Consequently, by the Banach--Steinhaus theorem, the norms of these operators are uniformly bounded; that is, there exist numbers $\varepsilon$ and $\gamma_{2}$ such that
$$
\left\|(L(\mu),B(\mu))^{-1}\right\|\leq\gamma_{2} \quad\mbox{for all} \  \mu\in\mathcal{B}(\mu_0, \varepsilon).
$$
Therefore, for all $\mu\in\mathcal{B}(\mu_0, \varepsilon)$ we obtain
\begin{gather*}
\big\|y(\cdot,\mu_0)-y(\cdot,\mu)\big\|_{(n+r)}
=\|(L(\mu),B(\mu))^{-1}(L(\mu),B(\mu))
(y(\cdot,\mu_0)-y(\cdot,\mu))\|_{(n+r)}\leq\\
\leq\gamma_{2}\,\|(L(\mu),B(\mu))
(y(\cdot,\mu_0)-y(\cdot,\mu))\|_{(C^{(n)})^{m}\times\mathbb{C}^{rm}}=\gamma_{2}\,\widetilde{d}_{(n)}(\mu).
\end{gather*}
This establishes the right-hand inequality in \eqref{3.6.bound}.
\end{proof}

The two-sided estimate \eqref{3.6.bound} characterizes the rate of convergence of the solutions $y(\cdot;\mu)$ as $\mu \rightarrow \mu_0$ in the sense of Ulam--Hyers stability (see, e.g., \cite{Rus,Hyers}).

\section{Limit theorems concerning the operators of the boundary-value problem}\label{Sec.3}

The main goal of this section is to investigate the continuity of the boundary operator $B(\mu)$, as well as the operator $\big(L(\mu),B(\mu)\big)$, with respect to the parameter in the uniform operator topology.
First, we establish limit theorems for the linear differential operator $L(\mu)$ as it tends to $L(\mu_0)$.

\begin{theorem}\label{eqival ym L}
The following assertions are equivalent as $\mu\to\mu_0$ in the metric space $\mathcal{M}$:
\begin{itemize}
\item [(I)] $A_{r-\ell}(\cdot,\mu) \rightarrow A_{r-\ell}(\cdot, \mu_0)$ in the space $(C^{(n)})^{m\times m}$ for each $\ell\in\{1,\ldots, r\}$;
\item [(II)] $L(\mu) \rightarrow L(\mu_0)$ in the uniform operator topology;
\item [(III)] $L(\mu) \rightarrow L(\mu_0)$ in the strong operator topology.
\end{itemize}
\end{theorem}

\begin{proof}[Proof of Theorem \ref{eqival ym L}]
It suffices to prove the implications: $(\mathrm{I})\Rightarrow(\mathrm{II})$ and $(\mathrm{III})\Rightarrow(\mathrm{I}).$
Suppose that $(\mathrm{I})$ holds. Then for every $y \in (C^{(n+r)})^m$, since $C^{(n)}$ is a Banach algebra, we have:
\begin{align*}
  \|(L(\mu) - L(\mu_0))y\|_{(n)}
  &\le \sum_{\ell=1}^r \|(A_{r-\ell}(\cdot, \mu) - A_{r-\ell}(\cdot, \mu_0)) y^{(r-\ell)}\|_{(n)} \\
  &\le c(n) \sum_{\ell=1}^r \|A_{r-\ell}(\cdot, \mu) - A_{r-\ell}(\cdot, \mu_0)\|_{(n)} \|y^{(r-\ell)}\|_{(n)} \\
  &\le c(n) \left(\sum_{\ell=1}^r \|A_{r-\ell}(\cdot, \mu) - A_{r-\ell}(\cdot, \mu_0)\|_{(n)}\right) \|y\|_{(n+r)}.
\end{align*}
This implies $(\mathrm{II})$.

Now suppose that $(\mathrm{III})$ holds. Then
\begin{equation}\label{vusok p}
\sum\limits_{l=1}^r A_{r-l}(\cdot,\mu)Y^{(r-l)}\to
\sum\limits_{l=1}^r A_{r-l}(\cdot, \mu_0)Y^{(r-l)},
\end{equation}
in $(C^{(n)})^{m\times m}$ for every $Y \in(C^{(n+r)})^{m\times m}$. Setting $Y(t):=I_m$, we obtain $A_{0}(\cdot,\mu)\to A_{0}(\cdot, \mu_0)$. Next, setting $Y(t):=tI_m$ in \eqref{vusok p}, we find that
\begin{equation*}
A_{1}(t,\mu)+A_{0}(t,\mu)t\to A_{1}(t, \mu_0)+A_{0}(t, \mu_0)t,
\end{equation*}
which implies that $A_{1}(\cdot,\mu)\to A_{1}(\cdot, \mu_0)$. Setting $Y(t):=t^{2}I_m$, we obtain
\begin{equation*}
2A_{2}(t,\mu)+2A_{1}(t,\mu)t+A_{0}(t,\mu)t^{2}\to 2A_{2}(t, \mu_0)+2A_{1}(t, \mu_0)t+A_{0}(t, \mu_0)t^{2}
\end{equation*}
and, consequently, $A_{2}(\cdot,\mu)\to A_{2}(\cdot, \mu_0)$. Continuing this process inductively, we arrive at $(\mathrm{I})$.
\end{proof}

We now examine the properties of the generic boundary operator $B(\mu)$. They will also be used in Section~\ref{Sec.4}.

\begin{theorem}\label{opBinSobolev}
Let a point $t_0 \in [a, b]$, matrices $\big(\alpha_{s}\big)_{s=1}^{n+r} \subset \mathbb{C}^{rm\times m}$, and a matrix-valued function $\Phi(\cdot)$ of dimension $rm \times m$ of bounded variation on $[a, b]$, continuous from the right on $(a, b)$ and vanishing at $t = a$, be given.
\begin{itemize}
\item [(1)] The operator
\begin{equation} \label{diyaBnash}
By=\sum _{s=1}^{n+r} \alpha_{s}\,y^{(s-1)}(t_0)+\int_{a}^b ({\rm d}\Phi(t))y^{(n+r)}(t), \quad y(\cdot)\in (C^{(n+r)})^{m},
\end{equation}
acts continuously from $(C^{(n+r)})^{m}$ to $\mathbb{C}^{rm}$, and its norm satisfies
\begin{equation*}\label{normaB}
\|B\|\leq \gamma\max  \limits_s \{\|\alpha_{s}\|\} + \| V_a^b \Phi(\cdot) \|_{\mathbb{C}^{rm \times m}},
\end{equation*}
where $\gamma > 0$ is a constant independent of $\alpha_{s}$ and $\Phi(\cdot)$.
\item [(2)] Every continuous linear operator $B\colon (C^{(n+r)})^m\rightarrow\mathbb{C}^{rm}$ admits a unique canonical representation of the form \eqref{diyaBnash}.
\end{itemize}
\end{theorem}

\begin{proof}[Proof of Theorem \ref{opBinSobolev}]
We first prove assertion $(1)$. Since the space $(C^{(n+r)})^{m}$ is continuously embedded in $(C^{(n+r-1)})^{m}$, there exists a constant $\gamma > 0$ such that
$$
\|y\|_{C^{(n+r-1)}}\leq \gamma \|y\|_{(n+r)}.
$$
Taking into account that $\Vert \cdot \Vert_\infty$ denotes the supremum norm, we have:
\begin{gather*}
\left\vert  \sum_{s=1}^{n+r}\alpha_s y^{(s-1)}(t_0) \right \vert
\leq  \sum_{s=1}^{n+r}\big\|\alpha_s \big\| \big|y^{(s-1)}(t_0)\big| \leq
\sum_{s=1}^{n+r}\big\|\alpha_s \big\| \big\|y^{(s-1)} \big\|_\infty
\leq \\ \max \limits_s \{\big\|\alpha_s \big\|\}\big\|y\big\|_{C^{(n+r-1)}}\leq \gamma \big\|y\big\|_{C^{(n+r)}}.
\end{gather*}
For the Riemann--Stieltjes integral, we obtain the estimate:
$$
\left| \int_{a}^b ({\rm d}\Phi(t))y^{(n+r)}(t) \right| \leq \| V_a^b \Phi(\cdot) \|_{\mathbb{C}^{rm \times m}} \left\|y^{(n+r)}\right\|_{\infty} \leq \| V_a^b \Phi(\cdot) \|_{\mathbb{C}^{rm \times m}} \left\|y\right\|_{C^{(n+r)}}.
$$
Therefore, the norm of $B$ satisfies
$$
\|B\|\leq \gamma \max \limits_s\{\big\|\alpha_s \big\|\} + \| V_a^b \Phi(\cdot) \|_{\mathbb{C}^{rm \times m}}.
$$

To prove assertion $(2)$, we recall the three-space factorization lemma (see \cite[Ch.~IV, Sect.~5]{KolmFom}): if $E$, $E_1$, and $E_2$ are linear spaces, and $B \colon E \rightarrow E_1$ and $A \colon E\rightarrow E_2$ are linear operators such that
$$
AE=E_2, \quad \ker B \supset \ker A,
$$
then there exists a unique linear operator $C \colon E_2 \rightarrow E_1$ such that $B=CA$.

Let $B\colon (C^{(n+r)})^m\rightarrow\mathbb{C}^{rm}$ be a continuous linear operator.
We denote by $(C^{(n+r)}_{0})^{m}$ the subspace of functions $f \in (C^{(n+r)})^{m}$ such that $f^{(s-1)}(t_0)=0$ for all $s\in \{1, \ldots, n+r\}$. We define the operator $A \colon (C^{(n+r)}_{0})^{m}\rightarrow (C^{(0)})^{m}$ by $A f:=f^{(n+r)}$.
Setting $E:=(C^{(n+r)}_{0})^m$, $E_1:=\mathbb{C}^{rm}$, and $E_2:=(C^{(0)})^{m}$,
let $\widetilde B:= B \mid (C^{(n+r)}_{0})^m$ denote the restriction of $B$ to $(C^{(n+r)}_{0})^m$. By the three-space lemma, there exists a unique linear operator $C \colon (C^{(0)})^{m}\rightarrow\mathbb{C}^{rm}$ such that $\widetilde Bg=CAg$ for all $g \in (C^{(n+r)}_{0})^m$.

By the classical Riesz representation theorem on $C^{(0)}$ (see, e.g., \cite[p.~344]{Dunford}), there exists a unique matrix-valued function $\Phi(\cdot)$ of dimension $rm \times m$ of bounded variation on $[a, b]$, continuous from the right on $(a, b)$ and vanishing at $t = a$, such that $C h = \int_a^b ({\rm d}\Phi(t)) h(t)$ for every $h \in (C^{(0)})^{m}$ (in the sense of the Riemann--Stieltjes integral). Hence,
\begin{equation*}
\widetilde B g= C Ag=\int_{a}^b ({\rm d}\Phi(t))g^{(n+r)}(t) \quad  \mbox{for all} \ g  \in (C^{(n+r)}_{0})^m.
\end{equation*}

Next, every $f \in (C^{(n+r)})^m$ can be represented as
$$
f(t)=\sum_{s=1}^{n+r}\frac{(t-t_0)^{s-1}}{(s-1)!}f^{(s-1)}(t_0)+g(t), \quad \mbox{where} \ g \in (C^{(n+r)}_{0})^m.
$$

Let $e_j$, $j=1,\ldots , m$, be the standard basis in $\mathbb{C}^{m}$, and set $h_{j,s} = B\big( \frac{(t-t_0)^{s-1}}{(s-1)!} e_j\big) \in \mathbb{C}^{rm}$ for each $s= 1, \ldots , n+r$. For each $s$, we define the matrix $\alpha_s$ of dimension $rm \times m$ by the condition:
\begin{equation*}
\alpha_s \Big( \sum_{j=1}^{m} x_{j,s}  e_j  \Big) = \sum_{j=1}^{m} x_{j,s}  h_{j,s}.
\end{equation*}
Choosing $x_{j,s}$ as the coordinates of $f^{(s-1)} (t_0)$ in the basis $\{e_j\}_{j=1,\ldots , m}$, we obtain
\begin{equation*}
B \left( \sum_{s=1}^{n+r}\frac{(t-t_0)^{s-1}}{(s-1)!}f^{(s-1)}(t_0) \right) =
\sum_{s=1}^{n+r}\alpha_s f^{(s-1)}(t_0)
\end{equation*}
and, consequently,
$$
Bf=\sum_{s=1}^{n+r}\alpha_s f^{(s-1)}(t_0)+\widetilde B(g) =
\sum_{s=1}^{n+r}\alpha_s f^{(s-1)}(t_0)+\int_{a}^b ({\rm d}\Phi(t))f^{(n+r)}(t).
$$

To prove uniqueness, suppose that $B$ also admits the representation
\begin{equation*}
By=\sum _{s=1}^{n+r} \beta_{s}\,y^{(s-1)}(t_0)+\int_{a}^b ({\rm d}\Psi(t))y^{(n+r)}(t), \quad y \in (C^{(n+r)})^{m} .
\end{equation*}
Then
\begin{equation*}
\sum _{s=1}^{n+r}\left(\alpha_s- \beta_{s}\right)\,y^{(s-1)}(t_0) =  \int_{a}^b ({\rm d}(\Psi(t)-\Phi(t)))y^{(n+r)}(t) \quad \mbox{for all} \  y \in (C^{(n+r)})^{m}.
\end{equation*}
This equality extends to matrix-valued functions:
\begin{equation}\label{JTnew1}
\sum _{s=1}^{n+r}\left(\alpha_s- \beta_{s}\right)\,Y^{(s-1)}(t_0)
= \int_{a}^b ({\rm d}(\Psi(t)-\Phi(t)))Y^{(n+r)}(t)
\end{equation}
for every $Y\in  (C^{(n+r)})^{m \times m}.$
The right-hand side of \eqref{JTnew1} vanishes if $Y$ is a matrix polynomial of degree at most $n+r-1$. For any $s , \ell \in \{1, \ldots, n+r \}$, there exists a matrix polynomial such that
$$
Y^{(\ell-1)}(t_0)=\delta_{\ell s}I_m,
$$
where $\delta_{\ell s}$ is the Kronecker delta.
Indeed, this follows from the Taylor formula at $t_0$: each entry $y_{i,j}$ of $Y=\left(y_{i,j}\right)_{i,j=1}^m$ is a scalar polynomial of degree $\leq n+r-1$ such that its derivatives $y_{i,j}^{(s-1)}(t_0)$ take prescribed values. Thus, $\alpha_s = \beta_{s}$ for all $s \in \{1, \ldots, n+r \}.$

It remains to show that $\Psi(t)=\Phi(t)$ for all $t \in [a,b]$. Using \eqref{JTnew1}, we deduce that
\begin{equation}\label{formyla6}
\int_{a}^b ({\rm d}(\Psi(t)-\Phi(t)))Z(t)=O_{rm \times m}
\end{equation}
for all $Z \in (C^{(0)})^{m \times m}.$
When $m=1$, this follows directly from the uniqueness theorem for linear continuous functionals on $C[a,b]$ (by F.~Riesz, normalized functions of bounded variation coinciding everywhere except possibly at points of discontinuity must be identically equal). For $m \geq 2$, we consider the $m^2$ matrix units $E_{i,j} \in \mathbb{C}^{m \times m}$, $i,j=1, \ldots,m$. Substituting $Z = \varphi E_{i,j}$ for arbitrary $\varphi \in C^{(0)}$ into \eqref{formyla6}, we conclude that each entry of the matrix-valued function $\Psi - \Phi$ is identically zero on $[a,b]$, which completes the proof.
\end{proof}

We now examine the continuity properties of the family of operators $B(\mu)$. Consider the following conditions as $\mu\to\mu_0$:

\begin{enumerate}
\renewcommand{\labelenumi}{\alph{enumi})}
\renewcommand{\theenumi}{\alph{enumi})}
\item[(a)] $\alpha_s(\mu)\rightarrow\alpha_s(\mu_0)$ in $\mathbb{C}^{rm\times m}$ for each
$s\in\{1, \ldots, n+r\}$;
\item[(b)]
$\| V_a^b \Phi(\cdot, \mu) \|_{\mathbb{C}^{rm \times m}} = O(1)$;
\item[(c)] $\Phi(b, \mu) \to \Phi(b, \mu_0)$;
\item[(d)] $\int_a^t \Phi(\tau,\mu)d\tau \rightarrow \int_a^t \Phi(\tau, \mu_0)d\tau$ for each $t \in (a, b]$;
\item[(e)] $V_a^b (\Phi(\cdot, \mu) - \Phi(\cdot, \mu_0)) \rightarrow0$.
\end{enumerate}

\begin{theorem}\label{strongunifB}
The family of operators $B(\mu)$ converges strongly to $B(\mu_0)$ as $\mu\to\mu_0$ if and only if conditions $(a)$, $(b)$, $(c)$, and $(d)$ hold. Moreover, $B(\mu)$ converges to $B(\mu_0)$ in the operator norm if and only if conditions $(a)$ and $(e)$ hold.
\end{theorem}

Clearly, the convergence is considered in $\mathbb{C}^{rm \times m}$. Condition $(e)$ is strictly stronger than conditions $(b)-(d)$. Indeed, $(e)$ implies the uniform convergence of $\Phi(t, \mu)$ to $\Phi(t, \mu_0)$ on $[a, b]$ as $\mu \to \mu_0$, whereas conditions $(b)-(d)$ do not even entail the pointwise convergence of these functions at any point of $(a, b)$.

\begin{proof}[Proof of Theorem \ref{strongunifB}] \textit{Sufficiency of $(a)$ and $(e)$ for uniform convergence.}
From the representation of $B(\mu)$ in Theorem~\ref{opBinSobolev}, we obtain
\begin{align*}
\|B(\mu) - B(\mu_0)\|
&\le \gamma\max_s |\alpha_s(\mu) - \alpha_s(\mu_0)|
  + \bigl\|V_a^b\bigl(\Phi(\cdot,\mu) - \Phi(\cdot,\mu_0)\bigr)\bigr\|_{\mathbb{C}^{rm\times m}},
\end{align*}
which immediately implies the sufficiency of conditions $(a)$ and $(e)$.
Furthermore, conditions $(a)$ and $(e)$ clearly imply conditions $(a)$--$(d)$.

\textit{Sufficiency of $(a)$, $(b)$, $(c)$, and $(d)$ for strong convergence.}
It suffices to show that
\begin{gather}
\|B(\mu)y - B(\mu_0)y\|
  \le \sum_{s=1}^{n+r} |\alpha_s(\mu) - \alpha_s(\mu_0)|\,
       \bigl|y^{(s-1)}(t_0)\bigr| + \nonumber \\
  \left|\int_a^b \bigl(d\Phi(t,\mu) - d\Phi(t,\mu_0)\bigr)y^{(n+r)}(t)\right| \label{formula2}
\end{gather}
tends to zero as $\mu\to\mu_0$.
Condition~$(a)$ guarantees that the discrete sum converges to zero.
We now analyze the integral term.

First, consider the case $m=1$. By the classical Helly--Bray theorem for the dual space of $C^{(0)}[a,b]$, the sequence of linear functionals
\[
  h \mapsto \int_a^b d\varphi(t,\mu)\, h(t),
  \qquad h\in C^{(0)}[a,b],
\]
generated by normalized functions of bounded variation $\varphi(\cdot,\mu)$, converges weakly$^*$ to the functional generated by $\varphi(\cdot,\mu_0)$ as $\mu\to\mu_0$ if and only if:
\begin{itemize}
  \item the variations $V_a^b\varphi(\cdot,\mu)$ are uniformly bounded (condition~(b));
  \item $\varphi(b,\mu)\to\varphi(b,\mu_0)$ (condition~(c));
  \item $\int_a^t\varphi(\tau,\mu)d\tau \to \int_a^t \varphi(\tau,\mu_0)d\tau$ for all $t\in(a,b]$ (condition~(d)).
\end{itemize}
Since $y^{(n+r)}\in C^{(0)}[a,b]$, conditions $(b)$, $(c)$, and $(d)$ imply that
\[
  \int_a^b \bigl(d\Phi(t,\mu) - d\Phi(t,\mu_0)\bigr)y^{(n+r)}(t) \to 0
  \quad\text{as }\mu\to\mu_0.
\]

In the case $m\ge 2$, formula \eqref{formula2} applies componentwise: for any matrix-valued function $Y\in (C^{(0)})^{m\times m}$,
\[
  \int_a^b \bigl(d\Phi(t,\mu) - d\Phi(t,\mu_0)\bigr)Y(t) \to 0
  \quad\text{as }\mu\to\mu_0.
\]
Choosing the matrix units $E_{i,j}\in\mathbb{C}^{m\times m}$ and setting $Y = \varphi E_{i,j}$ for arbitrary $\varphi\in C^{(0)}$, we obtain the result for each entry of $\Phi$.
Thus, conditions $(a)$, $(b)$, $(c)$, and $(d)$ are sufficient for strong convergence.

\textit{Necessity of $(a)$, $(b)$, $(c)$, and $(d)$ for strong convergence.}
Suppose that $B(\mu)y \to B(\mu_0)y$ for every $y\in (C^{(n+r)})^m$. In matrix form,
\begin{equation}\label{eq:3.6_C}
\begin{split}
  &\sum_{s=1}^{n+r} \alpha_s(\mu)\,Y^{(s-1)}(t_0)
  + \int_a^b d\Phi(t,\mu)\,Y^{(n+r)}(t)
  \\
  &\qquad\to\;
  \sum_{s=1}^{n+r} \alpha_s(\mu_0)\,Y^{(s-1)}(t_0)
  + \int_a^b d\Phi(t,\mu_0)\,Y^{(n+r)}(t)
\end{split}
\end{equation}
for all matrix-valued functions $Y\in (C^{(n+r)})^{m\times m}$. We will show that conditions $(a)$, $(b)$, $(c)$, and $(d)$ must hold.
\medskip

\textit{Condition $(a)$.}
Substitute $Y(t) := \frac{(t-t_0)^{j-1}}{(j-1)!} I_m$ for $j\in\{1,\ldots,n+r\}$.
Then $Y^{(s-1)}(t_0) = \delta_{s,j}\,I_m$ for each $s \in \{1,\ldots,n+r\}$ and $Y^{(n+r)}\equiv 0$.
Hence, the integral term in \eqref{eq:3.6_C} vanishes identically, leaving the single term at $s=j$:
\[
  \alpha_j(\mu) \to \alpha_j(\mu_0), \quad \mu\to\mu_0,
\]
for each $j \in \{1,\ldots,n+r\}$. This proves condition~(a).
\medskip

\textit{Conditions $(b)$, $(c)$, and $(d)$.}
First, let $m=1$.
Since condition~(a) has been established, relation \eqref{eq:3.6_C} implies:
\[
  \int_a^b d\Phi(t,\mu)\,y^{(n+r)}(t)
  \;\to\;
  \int_a^b d\Phi(t,\mu_0)\,y^{(n+r)}(t),
  \qquad y\in C^{(n+r)}[a,b].
\]
Since $n+r\ge 1$, the derivative $y^{(n+r)}$ ranges over the entire space $C^{(0)}[a,b]$ as $y$ ranges over $C^{(n+r)}[a,b]$. Thus, the linear functionals $u \mapsto \int_a^b d\Phi(t,\mu)\,u(t)$ converge weakly$^*$ on $C^{(0)}[a,b]$.
By the Helly--Bray theorem and the normalization of $\Phi(\cdot,\mu)$ ($\Phi(a,\mu)=0$, right-continuous on $(a,b)$), weak$^*$ convergence on $C^{(0)}$ is equivalent to the validity of conditions (b), (c), and~(d).
For $m\ge 2$, testing with $Y(t) := E_{i,j}\,\varphi(t)$ for $\varphi\in C^{(n+r)}$ and invoking condition~(a) establishes the convergence for each entry of $\Phi(\cdot,\mu)$.
This proves the necessity of $(b)$, $(c)$, and $(d)$.

\textit{Necessity of $(e)$ for uniform operator convergence.}
Suppose that $B(\mu)\to B(\mu_0)$ in the operator norm:
\begin{equation}\label{eq:3.7_C}
  \|B(\mu) - B(\mu_0)\| \to 0 \quad\text{as }\mu\to\mu_0.
\end{equation}
Then strong convergence holds, so conditions $(a)$, $(b)$, $(c)$, and $(d)$ are satisfied.
It remains to prove $(e)$:
\[
  V_a^b\bigl(\Phi(\cdot,\mu) - \Phi(\cdot,\mu_0)\bigr) \to 0
  \quad\text{as }\mu\to\mu_0.
\]
First, let $m=1$.
For any function $y\in C^{(n+r)}[a,b]$ with $y^{(s-1)}(t_0) = 0$ ($s = 1, \dots, n+r$) and $u := y^{(n+r)}\in C^{(0)}$, formula \eqref{eq:3.7_C} and Theorem~\ref{opBinSobolev} give
\[
  \left|\int_a^b \bigl(d\Phi(t,\mu) - d\Phi(t,\mu_0)\bigr)u(t)\right|
  \le \|B(\mu) - B(\mu_0)\| \cdot \|y\|_{C^{(n+r)}}
  \to 0.
\]
Since $\|y\|_{C^{(n+r)}} \le C\|u\|_\infty$, taking the supremum over all $u\in C^{(0)}$ with $\|u\|_\infty\le 1$ and recalling that the norm of the functional $u\mapsto\int_a^b d\psi(t)\,u(t)$ on $C^{(0)}[a,b]$ is exactly the total variation $V_a^b\psi$, we obtain
\[
  V_a^b\bigl(\Phi(\cdot,\mu) - \Phi(\cdot,\mu_0)\bigr) \to 0.
\]
For $m\ge 2$, the result follows by applying the scalar estimate to each entry via the matrix units $E_{i,j}$. This completes the proof of Theorem~\ref{strongunifB}.
\end{proof}

\section{Approximation by solutions of multipoint boundary-value problems} \label{Sec.4}

In this section, we apply the general limit theorems to approximate the solution of an arbitrary boundary-value problem by solutions of elementary multipoint problems with polynomial coefficients.

Consider the limit boundary-value problem
\begin{equation}\label{eq:4.1}
  (L_0 y_0)(t)
  := y_0^{(r)}(t) + \sum_{\ell=1}^{r} A_{r-\ell,0}(t)\,y_0^{(r-\ell)}(t)
  = f_0(t),\quad t\in(a,b),
\end{equation}
\begin{equation}\label{eq:4.2}
  B_0 y_0
  := \sum_{s=1}^{n+r} \alpha_{s,0}\,y_0^{(s-1)}(t_0)
     + \int_a^b \bigl(d\Phi_0(t)\bigr)\,y_0^{(n+r)}(t)
  = c_0,
\end{equation}
where the matrix-valued functions $A_{r-\ell,0}(\cdot)\in (C^{(n)})^{m\times m}$, the vector-valued function $f_0(\cdot)\in (C^{(n)})^m$, the vector $c_0\in\mathbb{C}^{rm}$, the matrices $\alpha_{s,0}\in\mathbb{C}^{rm\times m}$, and the matrix-valued function $\Phi_0(\cdot)$ of dimension $rm\times m$ has bounded variation on $[a,b]$, is continuous from the right on $(a,b)$, and satisfies $\Phi_0(a)=0$.
The integral in \eqref{eq:4.2} is understood as a Riemann--Stieltjes integral.
By Theorem~\ref{opBinSobolev}, every continuous linear boundary condition for system \eqref{eq:4.1} can be uniquely written in the form \eqref{eq:4.2}, where $t_0$ is an arbitrary fixed point in $[a,b]$.

Consider also the approximating family of problems:
\begin{equation}\label{eq:4.3}
  (L_k y_k)(t)
  := y_k^{(r)}(t) + \sum_{\ell=1}^{r} A_{r-\ell,k}(t)\,y_k^{(r-\ell)}(t)
  = f_k(t),\quad t\in(a,b),
\end{equation}
\begin{equation}\label{eq:4.4}
  B_k y_k
  := \sum_{s=1}^{n+r} \alpha_{s,k}\,y_k^{(s-1)}(t_0)
     + \sum_{j=1}^{N(k)} \beta_{j,k}\,y_k^{(n+r)}(t_{j,k})
  = c_k,\quad k=1,2,3,\ldots
\end{equation}
Here, the matrix-valued functions $A_{r-\ell,k}$ belong to a dense subset $F$ of $(C^{(n)})^{m\times m}$, and $G$ is a dense subset of $(C^{(n)})^m$, with $f_k\in G$ for all $k$. Furthermore, $\alpha_{s,k},\,\beta_{j,k}\in\mathbb{C}^{rm\times m}$, the points $t_{j,k}$ belong to a dense subset $P$ of $[a,b]$, and
\[
  f_k \to f_0,\qquad c_k \to c_0\quad\text{as }k\to\infty.
\]

We are interested in the conditions under which problems \eqref{eq:4.3}, \eqref{eq:4.4} are well-posed and their solutions satisfy
\[
  y_k \to y_0\quad\text{in }(C^{(n+r)})^m\text{ as }k\to\infty.
\]
To apply Theorem~\ref{nep v}, we interpret the index $k$ as a parameter in a metric space $\mathcal{M}$.
Setting $\mathcal{M}=\mathbb{Z}_+ := \mathbb{N} \cup \{0\}$, we equip $\mathcal{M}$ with the metric:
\[
  d(n,m) = d(m,n) =
  \begin{cases}
    \displaystyle\left|\frac{1}{n}-\frac{1}{m}\right|, & n\neq 0,\; m\neq 0,\\[6pt]
    \dfrac{1}{n}, & m=0,\; n\in\mathbb{N},\\[6pt]
    0, & n=m=0.
  \end{cases}
\]
Then $0$ is the unique limit point of $(\mathbb{Z}_+,d)$ and $d(0,n)\to 0\iff n\to\infty$.
We now formulate the main approximation theorem of this section.

\begin{theorem}\label{thm:4.1}
Let the homogeneous limit boundary-value problem $(L_0,B_0)y=0$ have only the trivial solution in $(C^{(n+r)})^m$ \textup{(}Condition~\textup{(0)}\textup{)}.
Then there exists a sequence of well-posed multipoint boundary-value problems of the form~\eqref{eq:4.3},~\eqref{eq:4.4} with polynomial coefficients and polynomial right-hand sides, in which the operators $L_k$ and $B_k$ are chosen independently of $f_0$ and $c_0$, such that
\[
  (L_k,B_k)^{-1} \xrightarrow{\,s\,} (L_0,B_0)^{-1}
  \quad\text{and}\quad
  y_k\to y_0\;\text{ in }(C^{(n+r)})^m
  \quad\text{as }k\to\infty.
\]
\end{theorem}

Before proving Theorem~\ref{thm:4.1}, we establish an auxiliary approximation lemma.

\begin{lemma}\label{lem:4.2}
Let $\Phi_0\colon[a,b]\to\mathbb{C}^{rm\times m}$ be a matrix-valued function of bounded variation, continuous from the right on $(a,b)$, with $\Phi_0(a)=0$.
Let $P\subset[a,b]$ be a dense subset.
Then there exists a sequence $\{\Phi_k\}_{k=1}^\infty$ of step functions \textup{(}with jumps only at points in $P$\textup{)} such that conditions \textup{(b)}, \textup{(c)}, and \textup{(d)} of Theorem~$3.3$ are satisfied, that is\textup{:}
\begin{enumerate}
  \item[\textup{(b)}]
    $V_a^b\Phi_k(\cdot) = O(1)$\textup{;}
  \item[\textup{(c)}]
    $\Phi_k(b)\to\Phi_0(b)$\textup{;}
  \item[\textup{(d)}]
    $\displaystyle\int_a^t\Phi_k(\tau)\,d\tau
     \to \int_a^t\Phi_0(\tau)\,d\tau$
    for all $t\in(a,b]$.
\end{enumerate}
\end{lemma}

\begin{proof}[Proof of Lemma \ref{lem:4.2}]
It suffices to consider the scalar case $m=1$; the matrix case then follows entrywise.
For each $k\in\mathbb{N}$, choose a partition
$a = t_0^k < t_1^k < \cdots < t_{N(k)}^k = b$ with $t_j^k \in P$ such that $\delta_k := \max_j(t_{j+1}^k - t_j^k)\to 0$ as $k\to\infty$.
We define the step functions
\begin{equation}\label{eq:Phi-step}
  \Phi_k(t) := \Phi_0\!\left(t_j^k\right),
  \quad t\in\bigl[t_j^k,\, t_{j+1}^k\bigr),\quad j=0,\ldots,N(k)-1,
\end{equation}
and set $\Phi_k(b) := \Phi_0(b)$.
Since $t_0^k = a$ and $\Phi_0(a) = 0$, on $[a, t_1^k)$ we have $\Phi_k(t) = \Phi_0(a) = 0$. Thus, $\Phi_k(a) = 0$ and $\Phi_k$ is continuous from the right on $[a,b)$ with jump points exclusively at $t_1^k, \ldots, t_{N(k)}^k \in P$.
We now verify the variation of $\Phi_k$. The jump of $\Phi_k$ at $t_j^k$ is $\Phi_0(t_j^k) - \Phi_0(t_{j-1}^k)$, so:
\[
  V_a^b\Phi_k
  = \sum_{j=0}^{N(k)-1}
    \bigl|\Phi_0\!\left(t_{j+1}^k\right)-\Phi_0\!\left(t_j^k\right)\bigr|
  \le V_a^b\Phi_0 < \infty,
\]
which proves condition~$(b)$.
Condition~$(c)$ holds trivially since $\Phi_k(b) = \Phi_0(b)$ for all $k$.
Since $\Phi_0\in BV[a,b]\subset L^\infty[a,b]\subset L^1[a,b]$, $\Phi_0$ is continuous almost everywhere with respect to Lebesgue measure (having at most countably many jump discontinuities).
The construction \eqref{eq:Phi-step} guarantees that $\Phi_k(t)\to\Phi_0(t)$ at every continuity point of $\Phi_0$. By the Lebesgue dominated convergence theorem, $\Phi_k \to \Phi_0$ in $L^1[a,b]$, which yields for every $t\in(a,b]$:
\[
  \int_a^t \Phi_k(\tau)\,d\tau
  \;\longrightarrow\;
  \int_a^t \Phi_0(\tau)\,d\tau
  \quad\text{as }k\to\infty.
\]
This proves condition~$(d)$ and completes the proof of Lemma~\ref{lem:4.2}.
\end{proof}

\begin{proof}[Proof of Theorem~\ref{thm:4.1}]
Since problem $(L_0,B_0)$ satisfies Condition~(0), by Theorem~\ref{nep v} it suffices to construct a sequence $\{L_k\}$ of differential operators and a sequence $\{B_k\}$ of boundary operators of the form \eqref{eq:4.4} such that $(L_k,B_k)\to(L_0,B_0)$ in the strong operator topology.

\textit{Step 1: Approximation of $L_0$.}
By the Weierstrass approximation theorem, polynomials are dense in $C^{(n)}[a,b]$. By Theorem~3.1, choosing matrix polynomials $A_{r-\ell,k}(\cdot)$ in \eqref{eq:4.3} that converge uniformly to $A_{r-\ell,0}(\cdot)$ in $C^{(n)}$, we obtain
\[
  L_k \to L_0\quad\text{in the uniform (and hence strong) operator topology.}
\]

\textit{Step 2: Approximation of $B_0$.}
By Theorem~3.2, the operator $B_0$ can be uniquely represented as
\begin{equation*}\label{eq:B0-repr}
  B_0 y
  = \sum_{s=1}^{n+r}\alpha_{s,0}\,y^{(s-1)}(t_0)
    + \int_a^b \bigl(d\Phi_0(t)\bigr)\,y^{(n+r)}(t),
  \quad y(\cdot)\in (C^{(n+r)})^m,
\end{equation*}
where $\Phi_0(\cdot)$ has bounded variation, is continuous from the right on $(a,b)$, and satisfies $\Phi_0(a)=0$.
By Lemma~\ref{lem:4.2}, there exists a sequence of step functions $\{\Phi_k\}$ with jumps in $P$ satisfying conditions~(b), (c), and~(d) of Theorem~3.3.
Setting $\alpha_{s,k}:=\alpha_{s,0}$, Theorem~3.3 guarantees that the operators
\begin{equation}\label{eq:4.5}
  B_k y
  = \sum_{s=1}^{n+r}\alpha_{s,k}\,y^{(s-1)}(t_0)
    + \int_a^b \bigl(d\Phi_k(t)\bigr)\,y^{(n+r)}(t),
  \quad y(\cdot)\in (C^{(n+r)})^m,
\end{equation}
converge \emph{strongly} to $B_0$: $B_k y\to B_0 y$ for every $y\in(C^{(n+r)})^m$.

\textit{Step 3: Reduction to multipoint conditions.}
We show that \eqref{eq:4.5} is equivalent to the multipoint condition \eqref{eq:4.4}.
First, let $m=1$.
Each step function $\Phi_k$ can be written as
\[
  \Phi_k(t) = \sum_{j=1}^{N(k)} c_{j,k}\,\mathbf{1}_{[t_{j,k},b]}(t),
  \quad t\in[a,b],
\]
where $c_{j,k}\in\mathbb{C}$, the points $t_{j,k}\in P$ are the jump points of $\Phi_k$ in $(a,b]$, and $\Delta\Phi_k(t_{j,k}) := \Phi_k(t_{j,k})-\Phi_k(t_{j,k}^-) = c_{j,k}$.
Since $y^{(n+r)}\in C^{(0)}[a,b]$, the Riemann--Stieltjes integral with respect to the step function $\Phi_k$ reduces to the finite sum:
\begin{equation*}\label{eq:stieltjes-step}
  \int_a^b y^{(n+r)}(t)\,d\Phi_k(t)
  = \sum_{j=1}^{N(k)} \Delta\Phi_k(t_{j,k})\cdot y^{(n+r)}(t_{j,k})
  = \sum_{j=1}^{N(k)} \beta_{j,k}\,y^{(n+r)}(t_{j,k}),
\end{equation*}
where $\beta_{j,k} := \Delta\Phi_k(t_{j,k}) = c_{j,k}\in\mathbb{C}^{rm\times m}$.
Thus, $B_k$ has precisely the multipoint form \eqref{eq:4.4}.
For $m\ge 2$, the same reasoning applies componentwise by decomposing $\Phi_k$ in the basis of matrix units $\{E_{i,j}\}_{i,j=1}^m$.

\textit{Step 4: Convergence of solutions.}
We have established that $L_k\to L_0$ and $B_k\to B_0$ strongly.
Applying Theorem~2.2 to the sequence of problems \eqref{eq:4.3}, \eqref{eq:4.4}, we conclude that
\[
  (L_k,B_k)^{-1}\xrightarrow{\,s\,}(L_0,B_0)^{-1}
  \quad\text{and}\quad
  y_k\to y_0\;\text{ in }(C^{(n+r)})^m \quad\text{as } k\to\infty.
\]
This completes the proof of Theorem~\ref{thm:4.1}.
\end{proof}

A natural question arises: under what conditions can the operator $(L_0, B_0)$ be approximated by multipoint problems of the form \eqref{eq:4.3}, \eqref{eq:4.4} in the \emph{uniform operator topology}:
\[
  \|(L_k, B_k)^{-1} - (L_0, B_0)^{-1}\| \to 0, \quad k \to \infty,
\]
instead of the strong operator convergence in Theorem~\ref{thm:4.1}? The following theorem gives an exhaustive answer.

\begin{theorem}\label{thm:4.3}
Let the conditions of Theorem~\ref{thm:4.1} hold. Then the uniform operator convergence
\begin{equation}\label{eq:4.7}
  \|(L_k, B_k)^{-1} - (L_0, B_0)^{-1}\| \to 0 \quad \text{as } k \to \infty
\end{equation}
is possible if and only if the matrix-valued function $\Phi_0$ is a pure jump function on $[a, b]$, that is, it admits the series representation:
\begin{equation}\label{eq:jump-repr}
  \Phi_0(t) = \sum_{j=1}^{\infty} c_j\,\mathbf{1}_{[t_j, b]}(t), \quad t \in [a, b],
\end{equation}
where $\{t_j\}_{j=1}^\infty \subset (a, b]$ is a countable set of jump points, $c_j \in \mathbb{C}^{rm \times m}$, and $\sum_{j=1}^{\infty} \|c_j\| < \infty$.
\end{theorem}

\begin{proof}
Since $(L_0, B_0)$ is invertible and $\|L_k - L_0\| \to 0$ (by Step~1 of Theorem~\ref{thm:4.1}), condition \eqref{eq:4.7} is equivalent to
\[
  \|B_k - B_0\| \to 0 \quad \text{as } k \to \infty.
\]
Because $\alpha_{s,k} = \alpha_{s,0}$, by the Riesz theorem the operator norm is given by the total variation:
\begin{equation}\label{eq:4.8}
  \|B_k - B_0\| = V_a^b\bigl(\Phi_k(\cdot) - \Phi_0(\cdot)\bigr) \to 0 \quad \text{as } k \to \infty,
\end{equation}
which is condition~(e) of Theorem~3.3.

\emph{Sufficiency.} Suppose that $\Phi_0$ is a pure jump function of the form \eqref{eq:jump-repr}. For simplicity, assume first that $\{t_j\} \subset P$. Since the series of variations $\sum_{j=1}^\infty \|c_j\|$ converges, for any $\varepsilon > 0$ there exists $N(\varepsilon) \in \mathbb{N}$ such that
\[
  \sum_{j=N(\varepsilon)+1}^{\infty} \|c_j\| < \varepsilon.
\]
For each $k \in \mathbb{N}$, we choose $\Phi_k$ as the partial sum of the first $N(k)$ jumps (where $N(k) \to \infty$ as $k \to \infty$):
\[
  \Phi_k(t) := \sum_{j=1}^{N(k)} c_j\,\mathbf{1}_{[t_j, b]}(t).
\]
Then the difference $\Phi_0(t) - \Phi_k(t) = \sum_{j=N(k)+1}^{\infty} c_j\,\mathbf{1}_{[t_j, b]}(t)$ is also a pure jump function, and its total variation is bounded by the tail:
\[
  V_a^b(\Phi_k - \Phi_0) = \sum_{j=N(k)+1}^{\infty} \|c_j\| \to 0 \quad \text{as } k \to \infty.
\]
If some points $t_j \notin P$, the density of $P$ in $[a, b]$ allows each $t_j$ to be approximated by points $\tau_{j,k} \in P$ such that the total variation of the difference remains bounded by $2\sum_{j>N(k)} \|c_j\| \to 0$. Thus, \eqref{eq:4.8} holds and $\|B_k - B_0\| \to 0$, which, combined with $\|L_k - L_0\| \to 0$ and the invertibility of $(L_0, B_0)$, yields \eqref{eq:4.7}.

\emph{Necessity.} Suppose that \eqref{eq:4.7} holds. Then by \eqref{eq:4.8}, $V_a^b(\Phi_k - \Phi_0) \to 0$ as $k \to \infty$, where each $\Phi_k$ is a step function with finitely many jumps.
By the Jordan--Lebesgue decomposition theorem, every function $\Phi_0 \in BV[a, b]$ with $\Phi_0(a) = 0$ can be uniquely decomposed into the sum of a pure jump function $\Phi_0^{\mathrm{jump}}$ and a continuous function of bounded variation $\Phi_0^{\mathrm{cont}}$:
\[
  \Phi_0 = \Phi_0^{\mathrm{jump}} + \Phi_0^{\mathrm{cont}}, \quad \Phi_0^{\mathrm{cont}}(a) = 0.
\]
Since each approximating function $\Phi_k$ is a step function, it has no continuous component, so $\Phi_k^{\mathrm{cont}} \equiv 0$. The total variation is additive with respect to mutually singular discrete and continuous measures:
\[
  V_a^b(\Phi_k - \Phi_0) = V_a^b\bigl(\Phi_k - \Phi_0^{\mathrm{jump}}\bigr) + V_a^b\bigl(\Phi_0^{\mathrm{cont}}\bigr).
\]
Consequently,
\[
  0 \le V_a^b\bigl(\Phi_0^{\mathrm{cont}}\bigr) \le V_a^b(\Phi_k - \Phi_0) \to 0 \quad \text{as } k \to \infty.
\]
Since the term $V_a^b\bigl(\Phi_0^{\mathrm{cont}}\bigr)$ is independent of $k$, we must have $V_a^b(\Phi_0^{\mathrm{cont}}) = 0$. Since $\Phi_0^{\mathrm{cont}}(a) = 0$, this implies $\Phi_0^{\mathrm{cont}}(t) \equiv 0$ for all $t \in [a, b]$. Therefore, $\Phi_0 = \Phi_0^{\mathrm{jump}}$ is a pure jump function, which proves \eqref{eq:jump-repr}.
\end{proof}

Note that the rate of convergence of the solutions is governed by Theorem~\ref{3.6.th-bound v}.

\begin{remark}\label{rem:4.4}
The result of Theorem~\ref{thm:4.3} reveals an essential difference between the spaces $C^{(n+r)}$ and the Sobolev spaces $W_p^{n+r}$ (cf. Theorem~4.3 in \cite{AMT2026}).
In the Sobolev case with $p = 1$, the dual space is $L^\infty[a, b]$, and the condition $\|B_k - B_0\| \to 0$ reduces to the approximation of $\Phi_0$ in the norm of $L^\infty$. Since step functions are dense in $L^\infty$ under the essential supremum norm, in $W_1^{n+r}$ uniform approximation is possible for \emph{arbitrary} $\Phi_0$.
By contrast, in $C^{(n+r)}$ the dual space is the space of functions of bounded variation equipped with the total variation norm $V_a^b$. Since the subspace of step functions is not dense in $BV[a, b]$ with respect to total variation, in $C^{(n+r)}$ uniform approximation is \emph{impossible} for any boundary condition containing a non-trivial continuous component (such as integral conditions). For such conditions in $C^{(n+r)}$, only \emph{strong} approximation is achievable (Theorem~\ref{thm:4.1}).
\end{remark}

\begin{example}\label{ex:1}
\upshape
Consider the Sturm--Liouville equation on $(0, 1)$ with boundary conditions containing integer and/or fractional derivatives:
\begin{equation}\label{ex1}
  (Ly)(t) := -y''(t) + q(t)y(t) = f(t), \quad t \in (0, 1),
\end{equation}
\begin{equation}\label{ex2}
  By := \sum_{j=1}^{N} \int_{0}^{1} \beta_{j}(t)\,\bigl(\textmd{D}_{0+}^{(\ell_j)}y\bigr)(t)\,{\rm d}t = c \in \mathbb{C}^2,
\end{equation}
where $q(\cdot) \in C^{(1)}[0, 1]$, $f(\cdot) \in C^{(1)}[0, 1]$, the vector-valued functions $\beta_{j}(\cdot) \in (C[0, 1])^2$, and $\textmd{D}_{0+}^{(\ell_j)}y$ denotes the Caputo fractional derivative of order
\begin{equation*}
  0 \le \ell_1 < \ell_2 < \ldots < \ell_N \le 3.
\end{equation*}
Recall (see, e.g., \cite{Caputo}) that for $\ell \in (k-1, k]$, $k \in \{1, 2, 3\}$, the Caputo fractional derivative is defined by
\[
  \bigl(\textmd{D}_{0+}^{(\ell)}y\bigr)(t) := \frac{1}{\Gamma(k - \ell)} \int_{0}^{t} (t - \tau)^{k - \ell - 1} y^{(k)}(\tau)\,{\rm d}\tau \quad \text{for } \ell \notin \mathbb{N}_0,
\]
and $\textmd{D}_{0+}^{(\ell)}y := y^{(\ell)}$ for $\ell \in \{0, 1, 2, 3\}$. If $y \in C^{(3)}[0, 1]$, then $\textmd{D}_{0+}^{(\ell_j)}y \in C[0, 1]$, so the operator $B : C^{(3)}[0, 1] \to \mathbb{C}^2$ is bounded.
Here $m = 1$, the order of the differential equation is $r = 2$, and the smoothness parameter is $n = 1$; hence, the solution space is $C^{(n+r)}[0, 1] = C^{(3)}[0, 1]$.
Assume that the characteristic matrix $M(L, B) \in \mathbb{C}^{2 \times 2}$ of problem \eqref{ex1}, \eqref{ex2} is non-singular (Condition~(0) holds). Then by \cite[Theorem~2.4]{NovAM2023}, there exists a unique solution $y \in C^{(3)}[0, 1]$.

By Theorem~\ref{thm:4.1}, we can approximate problem \eqref{ex1}, \eqref{ex2} by systems with polynomial coefficients $q_k(t)$, polynomial right-hand sides $f_k(t)$, and multipoint boundary conditions of the form \eqref{eq:4.4}:
\begin{equation}\label{ex:Bk}
  B_k y := \alpha_{1,k} y(0) + \alpha_{2,k} y'(0) + \alpha_{3,k} y''(0) + \sum_{j=1}^{N(k)} \beta_{j,k}\,y'''(t_{j,k}) = c_k \in \mathbb{C}^2,
\end{equation}
where $\alpha_{s,k}, \beta_{j,k} \in \mathbb{C}^2$, and $\{t_{j,k}\} \subset \mathbb{Q} \cap [0, 1]$ are rational nodes, such that the solutions $y_k \in C^{(3)}[0, 1]$ to problem \eqref{ex1}, \eqref{ex:Bk} satisfy
\begin{equation*}
  \|y - y_k\|_{(3)} \to 0 \quad \text{as } k \to \infty,
\end{equation*}
provided $\|q - q_k\|_{(1)} \to 0$, $\|f - f_k\|_{(1)} \to 0$, and $c_k \to c$. In terms of the solution operators:
\begin{equation*}
  (L_k, B_k)^{-1} \xrightarrow{\,s\,} (L, B)^{-1} \quad \text{as } k \to \infty.
\end{equation*}
However, because the boundary condition \eqref{ex2} is integral with continuous weight functions $\beta_j(t)$, the corresponding matrix-valued function $\Phi_0$ is absolutely continuous (its jump component is identically zero). Hence, by Theorem~\ref{thm:4.3}, uniform operator approximation is impossible:
\begin{equation*}
  \|(L_k, B_k)^{-1} - (L, B)^{-1}\| \not\to 0 \quad \text{as } k \to \infty.
\end{equation*}
\end{example}

\section{Acknowledgment}

The work of the first named author was funded by the Academy of Finland, grant no. 367369.

The second named author would like to thank the Isaac Newton Institute for Mathematical Sciences, Cambridge, for support and hospitality during the ``Solidarity Program'' where work on this paper was undertaken. This work was supported by ``EPSRC grant no EP/R014604/1'' and by the Ministry of Education and Sciences of Ukraine, grant no. 0126U000898. The author wishes to thank the Department of Mathematics, King's College London, for their hospitality.


\begin{thebibliography}{33}

\bibitem{Adomian}
G. Adomian, Solving Frontier Problems of Physics: The Decomposition Method. Kluwer Academic Publishers, The Netherlands, 1994. https://doi.org/10.1007/978-94-015-8289-6.

\bibitem{Ashordia97}
M. Ashordia, Boundary Value Problems for Systems of Linear Ordinary Differential Equations, Doctoral thesis, Masaryk University, Brno, 1997.

\bibitem{Ashordia98}
M. Ashordia, Conditions of the existence and uniqueness of solutions of the multipoint boundary value problem for a system of generalized ordinary differential equations, Georgian Math. J. 5 (1998) 1--24. https://doi.org/10.1515/GMJ.1998.1.

\bibitem{Ashordia}
M. Ashordia, Criteria of correctness of linear boundary value problems for systems of generalized ordinary differential equations, Czechoslov. Math. J. 46 (1996) 385--404.

\bibitem{Ashordia05}
M. Ashordia, On the general and multipoint boundary value problems for linear systems of generalized ordinary differential equations, linear impulse and linear difference systems, Mem. Differ. Equ. Math. Phys. 36 (2005) 1--80.

\bibitem{Ashordia20}
M. Ashordia, The general boundary value problems for linear systems of generalized ordinary differential equations, linear impulsive differential and ordinary differential systems. Numerical solvability, Mem. Differ. Equ. Math. Phys. 81 (2020) 1--184.


\bibitem{Atlasiuk2021}
O.M. Atlasiuk, Limit theorems for solutions of multipoint boundary-value problems with a parameter in Sobolev spaces, Ukrain. Math. J. 72 (2021) 1175--1184. https://doi.org/10.1007/s11253-020-01859-x.

\bibitem{AMT2026}
O. Atlasiuk, V. Mikhailets, J. Taskinen, Parameter-dependent inhomogeneous boundary-value problems in Sobolev spaces, (2025) https://doi.org/10.48550/arXiv.2512.21361


\bibitem{AtlasiukMikhailets20192}
O.M. Atlasiuk, V.A. Mikhailets, Fredholm one-dimensional boundary-value problems with parameter in Sobolev spaces, Ukrain. Math. J. 70 (2019) 1677--1687. https://doi.org/10.1007/s11253-019-01599-7.




\bibitem{Boyko}
V.M. Boyko, M. Kunzinger, R.O. Popovych, Parameter-dependent linear ordinary differential equations and topology of domains, J. Differ. Equ. 284 (2021) 546--575. https://doi.org/10.1016/j.jde.2021.03.001.

\bibitem{BrzdekPopaRus2018}
J.~Brzd{\k{e}}k, D.~Popa, I.~A.~Rus, N.~Fabiano, Ulam Type Stability, Springer, Cham, 2018.

\bibitem{Cartan1971}
H. Cartan, Differential Calculus, Houghton Mifflin, Boston / Hermann, Paris (1971).

\bibitem{Chua}
S.-K. Chua, Average value problems in ordinary differential equations, J. Differ. Equ. 249 (2010) 1531--1548. https://doi.org/10.1016/j.jde.2010.07.011.

\bibitem{Dunford}
N. Dunford, J.T. Schwartz, Linear Operators. Part I: General Theory, Interscience Publishers, New York–London, 1958.

\bibitem{Gikhman}
I.I. Gikhman, Concerning a theorem of N.N.~Bogolyubov, Ukr. Mat. Zh. 4 (1952) 215--219. (in Russian).

\bibitem{HRYP}
E.V. Gnyp, Continuity of the solutions of one-dimensional boundary-value problems with respect to the parameter in the Slobodetskii spaces, Ukrain. Math. J. 68 (2016) 846--861. https://doi.org/10.1007/s11253-016-1261-y.

\bibitem{HRYPKodliukMikhailets2015}
E.V. Gnyp, T.I. Kodlyuk, V.A. Mikhailets, Fredholm boundary-value problems with parameter in Sobolev spaces, Ukrain. Math. J. 67 (2015) 658--667. https://doi.org/10.1007/s11253-015-1105-1.


\bibitem{Goodrich}
C. Goodrich, C. Lizama, Existence and monotonicity of nonlocal boundary value problems: The one-dimensional case, Proc. R. Soc. Edinb. Sect. A 152 (2022) 1--27. https://doi.org/10.1017/prm.2020.90.

\bibitem{HRYPMikhailetsMurach2017}
Y. Hnyp, V. Mikhailets, A. Murach, Parameter-dependent one-dimensional boundary-value problems in Sobolev spaces, Electron. J. Differ. Equ. 2017 (2017) 1--13.


\bibitem{Horyunov}
A.S. Horyunov, Convergence and approximation of the Sturm--Liouville operators with potentials-distributions, Ukrain. Math. J. 67 (2015) 680--689. https://doi.org/10.1007/s11253-015-1107-z.

\bibitem{Hyers}
D.H. Hyers, G. Isac, T.M. Rassias, Stability of Functional Equations in Several Variables, Birkh\"{a}user Boston, MA, New York, 1998. https://doi.org/10.1007/978-1-4612-1790-9

\bibitem{Hermander1985}
L. H\"ormander, The Analysis of Linear Partial Differential Operators. III: Pseudo-differential operators, Springer-Verlag, Berlin, Heidelberg, 2007. https://doi.org/10.1007/978-3-540-49938-1.

\bibitem{Iha}
F.T. Iha, On boundary functionals and operators with finite-dimensional null spaces, Pac. J. Math. 56 (1975) 517--524. https://doi.org/10.2140/pjm.1975.56.517.

\bibitem{Jung2011}
S.-M.~Jung, Hyers--Ulam--Rassias Stability of Functional Equations in Nonlinear Analysis, Springer Optim. Appl., Vol.~48, Springer, New York, 2011.

\bibitem{Kato}
T.~Kato, Perturbation Theory for Linear Operators, Springer-Verlag, Berlin, Heidelberg, 1966.


\bibitem{Puza}
I. Kiguradze, B. P\.{u}\v{z}a, Boundary Value Problems for Systems of Linear Functional Differential Equations, Masaryk University, Brno, 2003.

\bibitem{Kigyradze1987}
I.T. Kiguradze, Boundary-value problems for systems of ordinary differential equations, J. Sov. Math. 43 (1988) 2259--2339. https://doi.org/10.1007/BF01100360.

\bibitem{Kigyradze2003}
I.T. Kiguradze, On boundary-value problems for linear differential systems with singularities, Differ. Equ. 39 (2003) 212--225. https://doi.org/10.1023/A:1025152932174.

\bibitem{Kigyradze1975}
I.T. Kiguradze, Some Singular Boundary-Value Problems for Ordinary Differential Equations, Tbilisi University, Tbilisi, 1975. (in Russian)	

\bibitem{Caputo}
A.A. Kilbas, H.M. Srivastava, J.J. Trujillo, Theory and Applications of Fractional Differential Equations, Elsevier, Amsterdam, 2006.

\bibitem{KodliukMikhailetsReva2013}
T.I. Kodlyuk, V.A. Mikhailets, N.V. Reva, Limit theorems for one-dimensional boundary-value problems, Ukrain. Math. J. 65 (2013) 77--90. https://doi.org/10.1007/s11253-013-0766-x.

\bibitem{KodliukMikhailets2013}
T.I. Kodlyuk, V.A. Mikhailets, Solutions of one-dimensional boundary-value problems with a parameter in Sobolev spaces, J. Math. Sci. 190 (2013) 589--599. https://doi.org/10.1007/s10958-013-1272-2.

\bibitem{KolmFom}
A.N. Kolmogorov, S.V. Fomin, Elements of the Theory of Functions and Functional Analysis, Dover Publications Inc., New York, 1999.	

\bibitem{Krall68}
A.M. Krall, Differential operators and their adjoints under integral and multiple point boundary conditions, J. Differ. Equ. 4 (1968)  327--336. https://doi.org/10.1016/0022-0396(68)90019-3.

\bibitem{Krall}
A.M. Krall, The development of general differential and general differential-boundary systems, Rocky Mt. J. Math. 5 (1975) 493--542. https://doi.org/10.1216/RMJ-1975-5-4-493.

\bibitem{KrasnKrein}
M.A. Krasnoselskii, S.G. Krein, On the principle of averaging in nonlinear mechanics, Uspekhi Mat. Nauk 10 (1955) 147--153. (in Russian)


\bibitem{Krein82}
S.G. Krein, Linear Equations in Banach Spaces, Birkh\"{a}user: Basel, Switzerland, 1982. https://doi.org/10.1007/978-1-4684-8068-9.

\bibitem{KurzweilVorel}
J. Kurzweil, Z. Vorel, Continuous dependence of solutions of differential equations on a parameter,  Czech. Math. J. 7 (1957) 568--583 (in Russian).

\bibitem{Levin}
A.Yu. Levin, The limiting transition for nonsingular systems, Dokl. Akad. Nauk SSSR 176 (1967) 774--777.

\bibitem{Ma}
R. Ma, A survey on nonlocal boundary value problems, Appl. Math. E-Notes 7 (2007) 257--279.







\bibitem{AtlasiukMikhailetsSkor}
V.A. Mikhailets, O.M. Atlasiuk, T.B. Skorobohach, On the solvability of Fredholm boundary-value problems in fractional Sobolev spaces, Ukrainian Math. J., 75 (2023)  107--117. https://doi.org/10.1007/s11253-023-02188-5. 

\bibitem{NovAM2023}
V. Mikhailets, O. Atlasiuk, The solvability of inhomogeneous boundary-value problems in Sobolev spaces, Banach J. Math. Anal. 18(2), 12 (2024). https://doi.org/10.1007/s43037-023-00316-8.


\bibitem{MikhMurachSol2016}
V.A. Mikhailets, A.A. Murach, V.O. Soldatov, Continuity in a parameter of solutions to generic boundary-value problems, Electron. J. Qual. Theory Differ. Equ. 2016, 87 (2016) 1--16. https://doi.org/10.14232/ejqtde.2016.1.87.

\bibitem{Pelekhata18}
V.A. Mikhailets, O.B. Pelekhata, N.V. Reva, Limit theorems for the solutions of boundary-value problems, Ukrainian Math. J. 70 (2018)  243--251. https://doi.org/10.1007/s11253-018-1498-8.


\bibitem{Nguen}
T.K. Nguen, On the dependence of a solution to a linear system of differential equations on a parameter, Differ. Equ. 29, 6 (1993)  830--835.

\bibitem{Opial}
Z. Opial, Continuous parameter dependence in linear systems of differential equations, J. Differ. Equ. 3 (1967) 571--579.

\bibitem{Paukstaite}
G. Pauk$\check{\mbox{s}}$tait\.e, A. $\check{\mbox{S}}$tikonas, Green's matrices for first order differential systems with nonlocal conditions, Math. Model. Anal. 22 (2017) 213--227. https://doi.org/10.3846/13926292.2017.1291456.

\bibitem{Pelekhata19}
O.B. Pelekhata, N.V. Reva, Limit theorems for the solutions of linear boundary-value problems for systems of differential equations, Ukrain. Math. J. 71 (2019) 1061--1070. https://doi.org/10.1007/s11253-019-01698-5.

\bibitem{PopaRasa2011}
D.~Popa, I.~Ra{\c{s}}a, On the Hyers--Ulam stability of the linear differential equation, J. Math. Anal. Appl. 381 (2011) 530--537.

\bibitem{Reid}
W.T. Reid, Some limit theorems for ordinary differential systems, J. Differ. Equ. 3 (1967) 423--439.

\bibitem{Rus}
I.A. Rus, Ulam Stability of the Operatorial Equations. In: T. Rassias, J. Brzdek, (eds) Functional Equations in Mathematical Analysis. Springer Optimization and Its Applications, vol 52. (2011) Springer, New York, NY. https://doi.org/10.1007/978-1-4614-0055-4\_23.

\bibitem{Skubachevskii}
A.L. Skubachevskii, Nonclassical boundary-value problems, J. Math. Sci. 155 (2008), 199–334. https://doi.org/10.1007/s10958-008-9218-9.

\bibitem{Smogorzhevsky}
O. Smogorzhevsky, Les fonctions de Green des syst\`{e}mes diff\'{e}rentielles lin\'{e}aires dans un domaine \`{a} une seule dimension, Mat. Sb. 7 (1940) 179--196.

\bibitem{Sold2015}
V.O. Soldatov, On the continuity in a parameter for the solutions of boundary-value problems total with respect to the spaces $C^{(n+r)}[a, b]$, Ukrainian Math. J. 67 (2015) 785--794. https://doi.org/10.1007/s11253-015-1114-0.

\bibitem{Soldatov2025}
V. Soldatov, Solvability of linear boundary-value problems for ordinary differential systems in the space $C^n$, Ukrainian Math. J. 77 (2025) 360--368. 
https://doi.org/10.1007/s11253-025-0251-x.

\bibitem{Whyburn}
W.M. Whyburn, Differential equations with general boundary conditions, Bull. Am. Math. Soc. 48 (1942) 692--704. https://doi.org/10.1090/S0002-9904-1942-07760-3.

\end{thebibliography}
\end{document}